\documentclass[12pt,twoside,american,english,british]{article}
\usepackage[T1]{fontenc}
\usepackage[utf8]{inputenc}
\usepackage[a4paper]{geometry}
\usepackage{color}
\usepackage{dsfont}
\usepackage{amsmath}
\usepackage{amsthm}
\usepackage{amssymb}
\usepackage{setspace}
\usepackage{esint}

\makeatletter
\numberwithin{equation}{section}
\numberwithin{figure}{section}
\theoremstyle{plain}
\newtheorem{thm}{\protect\theoremname}[section]
\theoremstyle{definition}
\newtheorem{defn}[thm]{\protect\definitionname}
\theoremstyle{plain}
\newtheorem{lem}[thm]{\protect\lemmaname}
\theoremstyle{plain}
\newtheorem{prop}[thm]{\protect\propositionname}
\newcommand{\lyxaddress}[1]{
	\par {\raggedright #1
	\vspace{1.4em}
	\noindent\par}
}

\usepackage[T1]{fontenc}
\usepackage[utf8]{inputenc}
\usepackage[a4paper]{geometry}
\usepackage{color}
\usepackage{dsfont}
\usepackage{amsmath}
\usepackage{amsthm}
\usepackage{amssymb}
\usepackage{setspace}
\usepackage{esint}

\makeatletter
\numberwithin{equation}{section}
\numberwithin{figure}{section}
\theoremstyle{plain}

\@ifundefined{date}{}{\date{}}
\usepackage{babel}
\usepackage{latexsym}
\usepackage{amsthm}
\usepackage{esint}
\usepackage{eucal}
\usepackage{epstopdf}
\usepackage{graphicx}
\usepackage{bigints}
\theoremstyle{plain}

\newtheoremstyle{boldremark}
    {\dimexpr\topsep/2\relax} 
    {\dimexpr\topsep/2\relax} 
    {}          
    {}          
    {\bfseries} 
    {.}         
    {.5em}      
    {}          

\theoremstyle{boldremark}
\newtheorem{brem} [thm] {Remark} 

\usepackage{fancyhdr}
\usepackage{blindtext}
\usepackage{titlesec}

\titleformat{name=\chapter}[display]
{\normalfont\Huge\itshape}
{\titlerule[1pt]\vspace{-33pt}\filleft%
  \parbox[t]{6em}{%
    \raggedleft%
    \rule{\linewidth}{0.5ex}\newline%
    \chaptertitlename\ \thechapter%
  }%
}
{4pc}
{\normalfont\upshape\bfseries\Huge}
\allowdisplaybreaks

\makeatother

\usepackage{babel}
\addto\captionsbritish{\renewcommand{\definitionname}{Definition}}
\addto\captionsbritish{\renewcommand{\lemmaname}{Lemma}}
\addto\captionsbritish{\renewcommand{\theoremname}{Theorem}}
\addto\captionsenglish{\renewcommand{\definitionname}{Definition}}
\addto\captionsenglish{\renewcommand{\lemmaname}{Lemma}}
\addto\captionsenglish{\renewcommand{\theoremname}{Theorem}}
\providecommand{\definitionname}{Definition}
\providecommand{\lemmaname}{Lemma}
\providecommand{\theoremname}{Theorem}

\usepackage[colorlinks,pdfpagelabels,pdfstartview = FitH,bookmarksopen
= true,bookmarksnumbered = true,linkcolor = blue,plainpages =
false,hypertexnames = false,citecolor = blue,pagebackref=false,urlcolor=blue]{hyperref}

\makeatother

\usepackage{babel}
\addto\captionsamerican{\renewcommand{\definitionname}{Definition}}
\addto\captionsamerican{\renewcommand{\lemmaname}{Lemma}}
\addto\captionsamerican{\renewcommand{\propositionname}{Proposition}}
\addto\captionsamerican{\renewcommand{\theoremname}{Theorem}}
\addto\captionsbritish{\renewcommand{\definitionname}{Definition}}
\addto\captionsbritish{\renewcommand{\lemmaname}{Lemma}}
\addto\captionsbritish{\renewcommand{\propositionname}{Proposition}}
\addto\captionsbritish{\renewcommand{\theoremname}{Theorem}}
\addto\captionsenglish{\renewcommand{\definitionname}{Definition}}
\addto\captionsenglish{\renewcommand{\lemmaname}{Lemma}}
\addto\captionsenglish{\renewcommand{\propositionname}{Proposition}}
\addto\captionsenglish{\renewcommand{\theoremname}{Theorem}}
\providecommand{\definitionname}{Definition}
\providecommand{\lemmaname}{Lemma}
\providecommand{\propositionname}{Proposition}
\providecommand{\theoremname}{Theorem}

\begin{document}
\title{\textbf{Gradient bounds for a widely degenerate orthotropic parabolic
equation}}
\author{Pasquale Ambrosio\thanks{\textbf{Corresponding author: Pasquale Ambrosio}, Department of Mathematics,
Uppsala University, P.O. Box 480, 751 06 Uppsala, Sweden.\textit{
E-mail address}: pasquale.ambrosio@math.uu.se}}
\date{June 1, 2026}
\maketitle
\begin{abstract}
\begin{singlespace}
\noindent In this paper, we consider the following nonlinear parabolic
equation 
\[
\partial_{t}u\,=\,\sum_{i=1}^{n}\partial_{x_{i}}\left[(\vert u_{x_{i}}\vert-\delta_{i})_{+}^{p-1}\frac{u_{x_{i}}}{\vert u_{x_{i}}\vert}\right]\,\,\,\,\,\,\,\,\,\,\mathrm{in}\,\,\,\Omega\times I,
\]
where $\Omega$ is a bounded open subset of $\mathbb{R}^{n}$ for
$n\geq2$, $I\subset\mathbb{R}$ is a bounded open interval, $p\geq2$,
$\delta_{1},\ldots,\delta_{n}$ are non-negative numbers and $\left(\,\cdot\,\right)_{+}$
denotes the positive part. We prove that the local weak solutions
are locally Lipschitz continuous in the spatial variable. The main
novelty here is that the above equation combines an orthotropic structure
with a strongly degenerate behavior. We emphasize that our result
can be considered, on the one hand, as the parabolic counterpart of
the elliptic result established in \cite{BBLV-elliptic}, and on the
other hand as an extension to a significantly more degenerate framework
of the findings contained in \cite{BBLV}. \vspace{0.2cm}
\end{singlespace}
\end{abstract}

\selectlanguage{english}%
\begin{singlespace}

\section{Introduction}
\end{singlespace}

\selectlanguage{british}%
\begin{singlespace}
\noindent $\hspace*{1em}$Let $\Omega\subset\mathbb{R}^{n}$ be a
bounded open set and $I\subset\mathbb{R}$ a bounded open interval.
We are interested in the gradient regularity of the local weak solutions
to the following parabolic equation
\begin{equation}
\partial_{t}u\,=\,\sum_{i=1}^{n}\partial_{x_{i}}\left[(\vert u_{x_{i}}\vert-\delta_{i})_{+}^{p-1}\frac{u_{x_{i}}}{\vert u_{x_{i}}\vert}\right]\,\,\,\,\,\,\,\,\,\,\mathrm{in}\,\,\,\Omega\times I,\label{eq:equation}
\end{equation}
where $p\geq2$, $\delta_{1},\ldots,\delta_{n}$ are non-negative
numbers and $\left(\,\cdot\,\right)_{+}$ stands for the positive
part. Throughout the paper, we denote by $T_{0}<T_{1}$ the endpoints
of the time interval $I$.\\
$\hspace*{1em}$Evolutionary equations of the above form have been
studied since the 1960s, notably by the Soviet school; see, for instance,
the work \cite{Vishik} by Vishik. Equation (\ref{eq:equation}) with
all $\delta_{i}$ set to zero is also explicitly presented in the
monographs \cite{Lions}, \cite[Example 4.A, Chapter III]{Show} and
\cite[Example 30.8]{Zeid}, among others.\\
$\hspace*{1em}$At first glance, (\ref{eq:equation}) looks quite
similar to the parabolic $p$-Laplace equation
\begin{equation}
\partial_{t}u\,=\,\sum_{i=1}^{n}\,(\vert Du\vert^{p-2}\,u_{x_{i}})_{x_{i}}\,\,\,\,\,\,\,\,\,\,\mathrm{in}\,\,\,\Omega\times I\,.\label{eq:parPLAP}
\end{equation}
However, the main novelty of equation (\ref{eq:equation}) lies in
the combination of two features, namely an \textit{orthotropic structure}
and a \textit{strongly degenerate behavior}. Indeed, unlike the parabolic
$p$-Laplace equation, for which the loss of ellipticity of the operator
$\mathrm{div}(\vert Du\vert^{p-2}Du)$ is restricted to a single point,
equation (\ref{eq:equation}) becomes degenerate on the larger set
\[
\bigcup_{i=1}^{n}\,\{\vert u_{x_{i}}\vert\leq\delta_{i}\}\,.
\]
$\hspace*{1em}$A more recent work in which equation (\ref{eq:equation})
appears with all $\delta_{i}$ equal to zero is \cite{BBLV}. There,
the authors derive local $L^{\infty}$ bounds for the spatial gradient
$Du$ of local weak solutions to (\ref{eq:equation}), but confining
their analysis to the case $p\geq2$ and $\max\,\{\delta_{i}\}=0$.
In this special case, as already observed in \cite{BBLV}, the basic
regularity theory equally applies to both (\ref{eq:equation}) and
(\ref{eq:parPLAP}). A classical reference in the field is DiBenedetto's
monograph \cite{DiBe}, which provides boundedness results for the
solution $u$ (see \cite[Chapter V]{DiBe}), Hölder continuity estimates
for $u$ (see \cite[Chapter III]{DiBe}), as well as Harnack inequalities
for non-negative solutions (see \cite[Chapter VI]{DiBe}). From a
technical point of view, there is no distinction to be made between
(\ref{eq:parPLAP}) and (\ref{eq:equation}) with all $\delta_{i}$
set to zero. Consequently, the results in \cite{BBLV} and \cite[Chapter V]{DiBe}
imply that, for $p\geq2$ and $\max\,\{\delta_{i}\}=0$, the local
weak solutions of (\ref{eq:equation}) are locally Lipschitz continuous
in the spatial variable.\\
$\hspace*{1em}$Concerning the gradient regularity of weak solutions
to equation (\ref{eq:parPLAP}), we refer again to DiBenedetto's book
for a comprehensive account of results on the subject, specifically
to \cite[Chapter VIII]{DiBe}. Since then, the literature on the regularity
for nonlinear, possibly degenerate or singular, parabolic equations
(or systems) has been steadily expanding, with the evolutionary $p$-Laplace
equation (\ref{eq:parPLAP}) serving as a prototypical model. Without
any claim of exhaustiveness, we can mention a few classical references
\cite{ChenDiBe,Choe,DiBe2,DiBeFri,DiBeFri2,Wieg}, up to more recent
contributions on the subject, such as \cite{BoDuMing,KuuMin1,KuuMin2},
among others.\\
$\hspace*{1em}$However, none of these results apply to equation (\ref{eq:equation}),
as they all rely on the fact that the loss of ellipticity of the operator
in divergence form is restricted to a single point, as in the model
case (\ref{eq:parPLAP}). As previously noted, such a property dramatically
fails for our equation (\ref{eq:equation}). Therefore, the aforementioned
references do not provide any regularity results for the spatial gradient
$Du$ of its solutions.\\
$\hspace*{1em}$In \cite{AmbCiani}, we have recently proved that
local weak solutions of (\ref{eq:equation}) are locally bounded also
in the case $p\geq2$ and $\max\,\{\delta_{i}\}>0$, thus extending
DiBenedetto's result \cite[Chapter V, Theorem 4.1]{DiBe} to our anisotropic
and more degenerate setting.\\
$\hspace*{1em}$The primary goal of this paper is to establish a local
$L^{\infty}$ bound on $Du$ for our equation (\ref{eq:equation}),
in order to extend the result of \cite{BBLV} to the markedly more
degenerate case $\max\,\{\delta_{i}\}>0$. To this end, we will need
to adapt the techniques developed in \cite{BB,BB2,BBJ,BBLV-elliptic,BBLV,BraCar,BraLePiVe}
for degenerate equations with orthotropic structure (see also \cite{Russo},
which deals with the higher differentiability of minimizers for non-autonomous
orthotropic functionals). Indeed, the main result of this work is
the following theorem, which can be considered as the parabolic counterpart
of the elliptic result \cite[Theorem 1.1]{BBLV-elliptic}. For notation
and definitions we refer to Section \ref{sec:prelim}.\medskip{}

\end{singlespace}
\begin{thm}
\label{thm:main}Let $n\geq2$ and $p\geq2$. Moreover, assume that
$u\in L_{loc}^{p}(I;W_{loc}^{1,p}(\Omega))$ is a local weak solution
of equation $(\ref{eq:equation})$. Then 
\[
Du\,\in\,L_{loc}^{\infty}(\Omega\times I,\mathbb{R}^{n})\,.
\]
More precisely, there exists a constant $C>1$, depending only on
$n$, $p$ and $\max\,\{\delta_{1},\ldots,\delta_{n}\}$, such that
for every parabolic cylinder $Q_{r}(x_{0},t_{0})\subset Q_{R}(x_{0},t_{0})\Subset\Omega\times I$
with $R\in(0,1]$, we have 
\begin{equation}
\Vert Du\Vert_{L^{\infty}(Q_{r}(x_{0},t_{0}))}\,\leq\,\frac{C}{(R-r)^{\vartheta p}}\left[1+\left(\iint_{Q_{R}(x_{0},t_{0})}\vert Du\vert^{p}\,dx\,dt\right)^{\frac{1}{2}}\right],\label{eq:main_est}
\end{equation}
where 
\[
\vartheta=\begin{cases}
\begin{array}{cc}
{\displaystyle \frac{n+2}{2}} & \,\,\mathit{if}\,\,n\geq3,\vspace{1mm}\\
\mathrm{\mathit{any}\,\,\mathit{number}}>2 & \,\,\mathit{if}\,\,n=2.
\end{array}\end{cases}
\]
\end{thm}

\noindent \smallskip{}
\begin{brem}We explicitly note that for $p\geq2$, the local weak
solutions of (\ref{eq:equation}) are locally Lipschitz continuous
in the spatial variable. Indeed, Theorem \ref{thm:main} provides
a local $L^{\infty}$ estimate for the spatial gradient $Du$. In
order to interpret this estimate as a genuine Lipschitz continuity
property for $u$, one also needs the local boundedness of the solution
itself. This is ensured by \cite[Theorem 1.1]{AmbCiani}, which holds
under the same assumptions on $p$ and $\delta_{1},\ldots,\delta_{n}$
considered here. We also observe that, as recalled in (\ref{eq:additional_properties-1}),
any local weak solution $u\in L_{loc}^{p}(I;W_{loc}^{1,p}(\Omega))$
of (\ref{eq:equation}) belongs to $C^{0}(\overline{J};L^{2}(\Omega'))$
for every subinterval $J\Subset I$ and every open set $\Omega'\Subset\Omega$,
so that the result of \cite{AmbCiani} is applicable in our framework.
Combining these facts, we obtain the claimed local Lipschitz continuity
with respect to the spatial variable.\end{brem}

\noindent \begin{brem}[\textbf{Comparison with other results}]More
generally, one may consider the following evolutionary PDE 
\begin{equation}
\partial_{t}u\,=\,\sum_{i=1}^{n}\partial_{x_{i}}\left[(\vert u_{x_{i}}\vert-\delta_{i})_{+}^{p_{i}-1}\frac{u_{x_{i}}}{\vert u_{x_{i}}\vert}\right]\,\,\,\,\,\,\,\,\,\,\mathrm{in}\,\,\,\Omega\times I,\label{eq:fully-aniso}
\end{equation}
which still exhibits an orthotropic structure. Now we have a full
range of exponents $1<p_{1}\leq p_{2}\leq\cdots\leq p_{n}$, one for
each spatial direction. We refer to \cite{Ters}, where some \textit{global}
Lipschitz regularity results are established for solutions to the
Cauchy-Dirichlet problem associated with (\ref{eq:fully-aniso}),
but only in the case $\max\,\{\delta_{i}\}=0$ and under suitable
regularity assumptions on the data. We emphasize that, due to their
global character, for $p_{1}=\cdots=p_{n}=p\geq2$, such results are
not comparable with the one proved here. We also mention \cite{CiaMosVes}
for a refined Harnack inequality for non-negative local weak solutions,
as well as for additional references on the problem. Furthermore,
in \cite{FeoVazVol} the authors provide an extensive analysis of
the Cauchy problem in the case $\max\,\{p_{i}\}<2$ and $\max\,\{\delta_{i}\}=0$,
along with related regularity results.

\noindent $\hspace*{1em}$However, as for the analog of our Theorem
\ref{eq:equation} for \textit{local solutions} of equation (\ref{eq:fully-aniso}),
this is still an open problem, as far as we know.\\
$\hspace*{1em}$It is worth noting that Theorem \ref{eq:equation}
can also be viewed as an extension to the orthotropic framework of
our previous result in \cite{AmbBau}, where we established local
$L^{\infty}$ bounds for the spatial gradient of weak solutions to
strongly degenerate parabolic systems, whose model case is given by
the equation
\begin{equation}
\partial_{t}u-\mathrm{div}\left((\vert Du\vert-\lambda)_{+}^{p-1}\frac{Du}{\vert Du\vert}\right)=f\,,\,\,\,\,\,\,\,\,\mathrm{with}\,\,\lambda>0\,.\label{eq:AmbPass}
\end{equation}
The main feature of this PDE is that the diffusion part is uniformly
elliptic only outside a ball with radius $\lambda$, while it behaves
asymptotically, that is, for large values of $\vert Du\vert$, like
the parabolic $p$-Laplace operator. Therefore, equations or systems
of the form (\ref{eq:AmbPass}) fall into the class of \textit{asymptotically
regular parabolic problems} (for a comprehensive overview of this
topic, see \cite{Amb1,AmbPass,BoDuGiPa} and the references therein).
As already pointed out in \cite{Amb2,BoDuGiPa}, no more than Lipschitz
regularity can be expected for solutions of equations or systems as
in (\ref{eq:AmbPass}). In fact, when $f=0$, any time-independent
$\lambda$-Lipschitz function solves (\ref{eq:AmbPass}), and even
more, it is a solution of the associated stationary equation or system.\\
$\hspace*{1em}$Finally, we mention the very recent contribution \cite{Strunk},
where the author studies parabolic equations of the type
\[
\partial_{t}u-\mathrm{div}\,D_{\xi}\mathcal{F}(x,t,Du)=g\,\,\,\,\,\,\,\,\,\,\mathrm{in}\,\,\,\Omega\times I,
\]
with $\mathcal{F}:\Omega\times I\times\mathbb{R}^{n}\to[0,\infty)$
satisfying the following conditions:\foreignlanguage{english}{\vspace{0.3cm}
}

\selectlanguage{english}%
\noindent $\hspace*{1em}(i)\,\,\,$$\mathcal{F}$ is only elliptic
for values of $Du$ outside a bounded and convex set $E\subset\mathbb{R}^{n}$
\foreignlanguage{british}{with the}\\
${\color{white}\hspace*{1em}(i)\,\,\,}$\foreignlanguage{british}{property
that $0\in\mathrm{Int}\,E$};\vspace{0.3cm}

\begin{singlespace}
\noindent $\hspace*{1em}(ii)\,\,\,$the partial map $\xi\mapsto\mathcal{F}(x,t,\xi)$
is regular on $\mathbb{R}^{n}\,\backslash\,\overline{E}$ and vanishes
whenever $\xi\in\overline{E}$.\vspace{0.3cm}

\end{singlespace}

\selectlanguage{british}%
\noindent Assuming that $g\in L^{n+2+\sigma}(\Omega\times I)$ for
some $\sigma>0$, in \cite{Strunk} it is shown that $\mathcal{K}(Du)\in C^{0}(\Omega\times I)$
for any function $\mathcal{K}\in C^{0}(\mathbb{R}^{n})$ that vanishes
on $E$. This result extends the $C^{1}$-regularity theorem proved
in \cite{BoDuGiPa} for equation (\ref{eq:AmbPass}) in the case $f=0$,
which fulfill $(i)$ and $(ii)$ with 
\[
\mathcal{F}(x,t,\xi)\,=\,\frac{1}{p}\,(\vert\xi\vert-\lambda)_{+}^{p}\,\,\,\,\,\,\,\,\,\,\,\,\mathrm{and}\,\,\,\,\,\,\,\,\,\,\,\,E=\{\xi\in\mathbb{R}^{n}:\vert\xi\vert<\lambda\}.
\]
\end{brem}

\subsection{Strategy of the proof\label{subsec:Strategy} }

\noindent $\hspace*{1em}$The key step in the proof of Theorem \ref{thm:main}
is to derive an \textit{a priori} Lipschitz estimate for smooth solutions
of the regularized parabolic equation
\begin{equation}
\partial_{t}u_{\varepsilon}\,=\,\mathrm{div}\,[D_{\xi}F_{\varepsilon}(Du_{\varepsilon})]\,,\label{eq:regul_equ}
\end{equation}
where $F_{\varepsilon}$ is a smooth, uniformly convex approximation
of the orthotropic function defined in (\ref{eq:F0}) below. Our goal
is to establish a local Lipschitz estimate for $u_{\varepsilon}$
that is independent of the regularization parameter $\varepsilon$.
Passing to the limit as $\varepsilon\to0$, we then show that the
family $\{u_{\varepsilon}\}$ converges to the original solution $u$.
This allows us to obtain the Lipschitz estimate for $u$ itself.\\
$\hspace*{1em}$In what follows, we emphasize the main difficulties
to get such a Lipschitz estimate on $u_{\varepsilon}$. At first sight,
the overall strategy may seem rather classical: we rely on a Moser
iterative scheme of reverse Hölder-type inequalities, arising from
the interplay between Caccioppoli estimates and Sobolev embeddings.
However, as precisely explained in \cite[Section 1.3]{BBLV}, this
standard approach cannot be applied directly, due to the severe degeneracy
of the orthotropic equation under consideration. For this reason,
we need to borrow the machinery developed in the elliptic setting
in \cite{BB} and \cite{BBLV-elliptic}, and later successfully exploited
in \cite{BBLV}.\\
$\hspace*{1em}$In particular, we make use of standard and nonstandard
(so-called ``weird'') Caccioppoli inequalities, both introduced
in \cite{BBLV} to deal with the solutions of the regularized equation
(\ref{eq:regul_equ}) (see Lemmas \ref{lem:standard} and \ref{lem:weird1}
below). Here, the term \textit{weird} refers to a family of parabolic
Caccioppoli-type estimates where two different components of the spatial
gradient $Du_{\varepsilon}$ are coupled. More precisely, such estimates
involve mixed terms depending simultaneously on $\partial_{x_{j}}u_{\varepsilon}$
and $\partial_{x_{k}}u_{\varepsilon}$, rather than a single spatial
derivative. This feature is crucial in the orthotropic setting, where
the degeneracy of the operator precludes the direct use of standard
elliptic estimates. By coupling different spatial directions, the
aforementioned “weird” Caccioppoli inequalities allow one to compensate
for the strong degeneracy of the orthotropic operator and recover
enough coercivity to carry out a Moser iterative scheme. Thanks to
these inequalities, we can proceed up to a certain point as in the
proof of \cite[Proposition 4.1]{BBLV}, where the parameters $\delta_{i}$
are all equal to $0$. In that setting, the combination of standard
and nonstandard Caccioppoli inequalities is sufficient to implement
the iteration scheme. A substantial difference arises in the present
framework, where $\max\,\{\delta_{i}\}>0$. In this case, the degeneracy
region of the orthotropic operator on the right-hand side of (\ref{eq:equation})
is no longer of measure zero, and this prevents a direct application
of the argument used in \cite{BBLV}. To overcome this difficulty
and handle the more degenerate situation where $\max\,\{\delta_{i}\}>0$,
we have to resort to an idea exploited in \cite[Section 5.2]{BBLV-elliptic},
in the elliptic setting. More specifically, we need to perform the
Moser iteration using a new measure $d\mu$, which is absolutely continuous
with respect to the $(n+1)$-dimensional Lebesgue measure and is supported
on the set whose complement coincides with the region\footnote{More precisely, in the limiting case $p=2$, the parameter $\delta$
is replaced by $\delta+1$.} 
\[
\left\{ \underset{1\,\le\,i\,\leq\,n}{\max}\left|\frac{\partial u_{\varepsilon}}{\partial x_{i}}\right|\leq\delta\right\} ,\,\,\,\,\,\,\,\,\,\,\mathrm{where}\,\,\,\delta=1+\max\,\{\delta_{1},\ldots,\delta_{n}\}\,.
\]
From the $L^{\infty}$ bound for $Du_{\varepsilon}$ obtained via
the Moser scheme with respect to the measure $d\mu$, we can easily
deduce the desired local $L^{\infty}$ estimate for $Du_{\varepsilon}$
in terms of the usual Lebesgue measure, by exploiting the support
properties of $d\mu$.\\
$\hspace*{1em}$From a technical viewpoint, the introduction of the
measure $d\mu$ allows the Moser iteration to be localized away from
the degeneracy set of the orthotropic operator. This modification
of the underlying measure in the iteration scheme makes the Moser
iterative procedure fully compatible with the orthotropic structure
of the equation and with the presence of multiple degeneracy thresholds
$\delta_{i}\geq0$. This idea represents the main new ingredient in
our approach compared to \cite{BBLV} and allows us to achieve the
proof of Theorem \ref{thm:main}.

\subsection{Plan of the paper }

\noindent $\hspace*{1em}$The paper is organized as follows. Section
\ref{sec:prelim} collects the preliminary material, including classical
notations, basic properties of local weak solutions to (\ref{eq:equation})
and a few auxiliary lemmas. In Section \ref{sec:Regularization},
we define the exact structure of the regularized equation (\ref{eq:regul_equ}),
whose local weak solution $u_{\varepsilon}$ satisfies the Caccioppoli-type
inequalities established in \cite{BBLV} and recalled in Subsection
\ref{subsec:Caccioppoli}. Section \ref{sec:Moser} is devoted to
the Moser iterative scheme, which yields a uniform $L^{\infty}$ bound
for $Du_{\varepsilon}$ through the uniform energy estimates derived
in Section \ref{sec:energy}. Finally, in Section \ref{sec:Dimostrazione}
we complete the proof of Theorem \ref{thm:main}, by transferring
the \textit{a priori} estimates obtained for the approximating solutions
$u_{\varepsilon}$ to the original solution $u$. In treating the
limiting case $p=2$, we also rely on a lemma proved in the \hyperref[sec:appendice]{Appendix}.
\selectlanguage{english}%
\begin{singlespace}

\section{Notation and preliminaries\label{sec:prelim}}
\end{singlespace}

\selectlanguage{british}%
\noindent $\hspace*{1em}$In this paper we shall denote by $C$ or
$c$ a general positive constant that may vary on different occasions,
even within the same line of estimates. Relevant dependencies on parameters
and special constants will be suitably emphasized using parentheses
or subscripts. The norm we use on $\mathbb{R}^{k}$, $k\in\mathbb{N}$,
will be the standard Euclidean one and it will be denoted by $\left|\,\cdot\,\right|$.
In particular, for the vectors $\xi,\eta\in\mathbb{R}^{k}$, we write
$\langle\xi,\eta\rangle$ for the usual inner product and $\left|\xi\right|:=\langle\xi,\xi\rangle^{\frac{1}{2}}$
for the corresponding Euclidean norm.\\
$\hspace*{1em}$In what follows, we use the notation 
\[
Q_{R}(x_{0}):=\,x_{0}+(-R,R)^{n},\,\,\,\,\,\,\,\,\,\,\,\,x_{0}\in\mathbb{R}^{n},\,R>0,
\]
for the $n$-dimensional open cube centered at $x_{0}$ with side
length $2R$. In the paper we will also work with \textit{anisotropic}
parabolic cylinders of the type
\begin{equation}
Q_{R}(x_{0},t_{0}):=\,Q_{R}(x_{0})\times(t_{0}-R^{p},t_{0}),\,\,\,\,\,\,\,\,\,\,\,\,t_{0}\in\mathbb{R}.\label{eq:anis_cylin}
\end{equation}

\selectlanguage{english}%
\noindent $\hspace*{1em}$If $E\subseteq\mathbb{R}^{k}$ is a Lebesgue-measurable
set, then we will denote by $\vert E\vert$ its $k$-dimensional Lebesgue
measure.\\
$\hspace*{1em}$In the sequel, we will also adopt the following notation
for convenience: given $v\in L^{1}(\Omega\times I)$, 
\[
\int_{\Omega\times\{\tau\}}v\,dx:=\int_{\Omega}v(x,\tau)\,dx\,\,\,\,\,\,\,\,\,\,\mathrm{for\,\,a.e.\,\,}\tau\in I.
\]
\foreignlanguage{british}{$\hspace*{1em}$}In this work, we define
a local weak solution to (\ref{eq:equation}) as follows.\medskip{}

\selectlanguage{british}%
\begin{defn}
\noindent Let $F_{0}:\mathbb{R}^{n}\to\mathbb{R}$ be the function
defined by 
\begin{equation}
F_{0}(\xi):=\sum_{i=1}^{n}\,\frac{1}{p}\,(\vert\xi_{i}\vert-\delta_{i})_{+}^{p}\,.\label{eq:F0}
\end{equation}
We say that $u\in L_{loc}^{p}(I;W_{loc}^{1,p}(\Omega))$ is a \textit{local
weak solution} of equation (\ref{eq:equation}) if and only if, for
any test function $\varphi\in C_{0}^{\infty}(\Omega\times I)$, we
have
\begin{equation}
\iint_{\Omega\times I}\left(u\,\partial_{t}\varphi\,-\langle D_{\xi}F_{0}(Du),D\varphi\rangle\right)dx\,dt\,=\,0\,.\label{eq:loc_weak_sol}
\end{equation}
\end{defn}

\noindent \smallskip{}
$\hspace*{1em}$Let $u\in L_{loc}^{p}(I;W_{loc}^{1,p}(\Omega))$ be
a local weak solution of (\ref{eq:equation}). We now recall some
additional properties of $u$: for every subinterval $J\Subset I$
and every open set $\Omega'\Subset\Omega$, we have 
\begin{equation}
\partial_{t}u\in L^{p'}(J;W^{-1,p'}(\Omega'))\,\,\,\,\,\,\,\,\,\,\,\,\,\,\,\mathrm{and}\,\,\,\,\,\,\,\,\,\,\,\,\,\,\,u\in C^{0}(\overline{J};L^{2}(\Omega'))\,.\label{eq:additional_properties-1}
\end{equation}
Here, $p':=p/(p-1)$ is the conjugate exponent of $p$ and $W^{-1,p'}(\Omega')$
is the topological dual space of $W_{0}^{1,p}(\Omega')$, which in
turn is the completion of $C_{0}^{\infty}(\Omega')$ with respect
to the $L^{p}$ norm of the gradient.\\
$\hspace*{1em}$In what follows, we briefly recall the argument used
to obtain (\ref{eq:additional_properties-1}). Fix $J$ and $\Omega'$
as above. By using equation (\ref{eq:loc_weak_sol}) and the fact
that $D_{\xi}F_{0}(Du)\in L_{loc}^{p'}(\Omega\times I,\mathbb{R}^{n})$,
one easily gets\begin{align*}
\left|\iint_{\Omega'\times J}u\,\partial_{t}\psi\,dx\,dt\right|&=\left|\iint_{\Omega'\times J}\langle D_{\xi}F_{0}(Du),D\psi\rangle\,dx\,dt\right|\\
&\leq\,n\,\Vert Du\Vert_{L^{p}(\Omega'\times J)}^{p-1}\,\Vert D\psi\Vert_{L^{p}(\Omega'\times J)}\\
&\leq\,n\,\Vert Du\Vert_{L^{p}(\Omega'\times J)}^{p-1}\,\Vert\psi\Vert_{L^{p}(J;W^{1,p}_{0}(\Omega'))},\,\,\,\,\,\,\,\,\,\,\mathrm{for\,\,every\,\,}\psi\in C_{0}^{\infty}(\Omega'\times J)\,.
\end{align*}By density, we can extend the linear functional 
\[
\Lambda:\psi\mapsto\iint_{\Omega'\times J}u\,\partial_{t}\psi\,dx\,dt
\]
to the whole space $L^{p}(J;W_{0}^{1,p}(\Omega'))$. This implies
that (see, for example, \cite[Theorem 1.5, Chapter III]{Show})
\[
\Lambda\in\left(L^{p}(J;W_{0}^{1,p}(\Omega'))\right)^{*}=\,L^{p'}(J;W^{-1,p'}(\Omega'))\,.
\]
From the definition of $\Lambda$ and of weak derivative, we obtain
the first property in (\ref{eq:additional_properties-1}).\\
$\hspace*{1em}$The second property in (\ref{eq:additional_properties-1})
follows by recalling that for every open set $\mathcal{O}\subset\mathbb{R}^{n}$,
we have (see \cite[Proposition 1.2, Chapter III]{Show}) 
\[
\mathcal{W}_{p}(\mathcal{O}\times J):=\left\{ v\in L^{p}(J;W_{0}^{1,p}(\mathcal{O}))\,:\,\partial_{t}v\in L^{p'}(J;W^{-1,p'}(\mathcal{O}))\right\} \subset C^{0}(\overline{J};L^{2}(\mathcal{O}))\,.
\]
Indeed, it is sufficient to take $\Omega'\Subset\Omega''\Subset\Omega$
and use the previous inclusion for the function $\eta u$, where $\eta\in C_{0}^{\infty}(\Omega'')$
is such that $\eta\equiv1$ on $\Omega'$. By the first property in
(\ref{eq:additional_properties-1}) (used with $\Omega''$ in place
of $\Omega'$) and the properties of $\eta$, we have that 
\[
\eta u\in\mathcal{W}_{p}(\Omega''\times J)\subset C^{0}(\overline{J};L^{2}(\Omega''))\,.
\]
Since $\eta\equiv1$ on $\Omega'$, the above fact implies that $u\in C^{0}(\overline{J};L^{2}(\Omega'))$,
as claimed.\bigskip{}

\noindent $\hspace*{1em}$We now gather some lemmas that will be useful
to prove our results. Let $J_{\lambda}:\mathbb{R}\to\mathbb{R}$ and
$H_{\lambda}:\mathbb{R}\to\mathbb{R}$ be the auxiliary functions
defined respectively by 
\begin{equation}
J_{\lambda}(s):=\begin{cases}
\begin{array}{cc}
{\displaystyle (\vert s\vert-\lambda)_{+}^{p-1}\,\frac{s}{\vert s\vert}} & \,\,\mathrm{if}\,\,\,s\neq0,\\
0 & \,\,\mathrm{if}\,\,\,s=0,
\end{array}\end{cases}\label{eq:J}
\end{equation}
and 
\begin{equation}
H_{\lambda}(s):=\begin{cases}
\begin{array}{cc}
{\displaystyle (\vert s\vert-\lambda)_{+}^{\frac{p}{2}}\,\frac{s}{\vert s\vert}} & \,\,\mathrm{if}\,\,\,s\neq0,\\
0 & \,\,\mathrm{if}\,\,\,s=0,
\end{array}\end{cases}\label{eq:H}
\end{equation}

\noindent where $\lambda\ge0$ is a parameter. We record the following
estimate, which can be obtained by suitably modifying the proof of
\cite[Lemma 4.1]{BraCaSan} in the case $N=1$.
\begin{lem}
\noindent \label{lem:Brasco}Let $p\geq2$ and $\lambda\geq0$. Then,
for every $s,t\in\mathbb{R}$ we get 
\[
\left(J_{\lambda}(s)-J_{\lambda}(t)\right)(s-t)\,\geq\,\frac{4}{p^{2}}\left|H_{\lambda}(s)-H_{\lambda}(t)\right|^{2}.
\]
\end{lem}

\noindent $\hspace*{1em}$The next result follows from an elementary
computation (see, e.g., \cite[Lemma A.1]{BraLePiVe}).
\begin{lem}
\noindent \label{lem:banal}Let $v:\mathbb{R}\rightarrow[0,\infty)$
be a $C^{1,1}$ convex function. Let us set 
\begin{equation}
\mathcal{G}(s)=\int_{0}^{s}\sqrt{v''(\tau)}\,d\tau\,.\label{eq:G-integrale}
\end{equation}
Then, for every $a,b\in\mathbb{R}$ we have 
\[
\left(v'(a)-v'(b)\right)(a-b)\,\geq\,\vert\mathcal{G}(a)-\mathcal{G}(b)\vert^{2}\,.
\]
\end{lem}

\noindent \begin{brem}Since the function $v$ in Lemma \ref{lem:banal}
belongs to $C^{1,1}(\mathbb{R})$, its derivative $v'$ is Lipschitz
continuous. Hence, $v''$ exists almost everywhere in $\mathbb{R}$
(by Rademacher's theorem) and belongs to $L^{\infty}(\mathbb{R})$.
In particular, the integral in (\ref{eq:G-integrale}) is well defined
as a Lebesgue integral.\end{brem}

\noindent $\hspace*{1em}$We conclude this section with the following
classical result; see \cite[Lemma 6.1]{Giu} for a proof. 
\begin{lem}
\begin{singlespace}
\noindent \label{lem:giusti} Let $0\leq\rho_{0}<\rho_{1}<\infty$
and assume that $Z:[\rho_{0},\rho_{1}]\rightarrow[0,\infty)$ is a
bounded function satisfying 
\[
Z(\rho)\,\leq\,\theta\,Z(r)\,+\,\frac{A}{(r-\rho)^{\alpha}}\,+\,\frac{B}{(r-\rho)^{\beta}}\,+\,C
\]
for all $\rho_{0}\leq\rho<r\leq\rho_{1}$, for some $\theta\in(0,1)$
and fixed non-negative constants $A$, $B$, $C$, $\alpha>\beta>0$.
Then, there exists a constant $\kappa=\kappa(\alpha,\theta)>0$ such
that 
\[
Z(\rho_{0})\,\leq\,\kappa\left(\frac{A}{(\rho_{1}-\rho_{0})^{\alpha}}\,+\,\frac{B}{(\rho_{1}-\rho_{0})^{\beta}}\,+\,C\right).
\]
\end{singlespace}
\end{lem}

\section{Estimates for a regularized equation\label{sec:Regularization} }

\noindent $\hspace*{1em}$We define 
\begin{equation}
G(\xi):=\,\frac{1}{p}\,(1+\vert\xi\vert^{2})^{\frac{p}{2}},\,\,\,\,\,\,\,\,\,\,\xi\in\mathbb{R}^{n},\label{eq:G}
\end{equation}
and, for $i\in\{1,\ldots,n\}$, we set 
\[
g_{i}(s):=\,\frac{1}{p}\,(\vert s\vert-\delta_{i})_{+}^{p}\,,\,\,\,\,\,\,\,\,\,\,s\in\mathbb{R}\,.
\]
For $p>2$ and for every $\varepsilon\in(0,1)$, we consider the convex
function
\begin{equation}
F_{\varepsilon}(\xi):=\sum_{i=1}^{n}g_{i}(\xi_{i})+\,\varepsilon\,G(\xi)\,,\,\,\,\,\,\,\,\,\,\,\xi\in\mathbb{R}^{n}.\label{eq:F_eps}
\end{equation}
\begin{brem}For $p=2$ and $\delta_{i}>0$, we have $g_{i}\in C^{1}(\mathbb{R})\cap C^{\infty}(\mathbb{R}\,\backslash\,\{-\delta_{i},\delta_{i}\})$,
but $g_{i}$ is not in $C^{2}(\mathbb{R})$. In this case, one would
need to replace $g_{i}$ by a regularized version $g_{i,\varepsilon}$,
in particular for the derivation of the ellipticity bounds in Lemma
\ref{lem:ellipticity} below. For $0<\varepsilon<\min\,\{1,\delta_{i}\}$,
a regularized version of $g_{i}$ is given by
\begin{equation}
g_{i,\varepsilon}(s):=\begin{cases}
\begin{array}{cc}
0 & \mathrm{if}\,\,\vert s\vert\le\delta_{i}-\varepsilon\,,\,\,\,\,\,\,\,\,\,\,\,\,\,\,\,\,\,\,\,\,\vspace{2mm}\\
{\displaystyle \frac{1}{12\,\varepsilon}\,(\vert s\vert-\delta_{i}+\varepsilon)^{3}} & \,\,\,\mathrm{if}\,\,\,\delta_{i}-\varepsilon\le\vert s\vert\le\delta_{i}+\varepsilon\,,\vspace{2mm}\\
{\displaystyle \frac{\varepsilon^{2}}{6}\,+\,\frac{1}{2}\,(\vert s\vert-\delta_{i})^{2}} & \mathrm{if}\,\,\vert s\vert\ge\delta_{i}+\varepsilon\,,\,\,\,\,\,\,\,\,\,\,\,\,\,\,\,\,\,\,\,\,
\end{array}\end{cases}\label{eq:g_eps}
\end{equation}
which converges in $C^{1}$ to $(\vert s\vert-\delta_{i})_{+}^{2}/2$
as $\varepsilon$ goes to $0$ (see Lemma \ref{lem:lemmaA1} in the
\hyperref[sec:appendice]{Appendix} and \cite[Section 2]{BraCar}).
Therefore, for $p=2$, to ensure that $F_{\varepsilon}\in C^{2}(\mathbb{R}^{n})$,
we need to define it as follows:\foreignlanguage{english}{\vspace{0.3cm}
}

\selectlanguage{english}%
\noindent $\hspace*{1em}\bullet\,\,\,$first, we denote $\Delta^{+}:=\{\delta_{i}:i\in\{1,\ldots,n\},\delta_{i}>0\}$;\vspace{0.3cm}

\begin{singlespace}
\noindent $\hspace*{1em}\bullet\,\,\,$then, for $\varepsilon\in(0,\min\,\{1,\inf\Delta^{+}\})$,
we define\vspace{0.3cm}

\end{singlespace}

\selectlanguage{british}%
\noindent 
\begin{equation}
F_{\varepsilon}(\xi)\,=\,\sum_{i=1}^{n}\,\widetilde{g}_{i,\varepsilon}(\xi_{i})\,+\,\frac{\varepsilon}{2}\,(1+\vert\xi\vert^{2})\,,\label{eq:F_quad}
\end{equation}
with the convention that $\inf\Delta^{+}=+\infty$ if $\Delta^{+}=\emptyset$
and 
\begin{equation}
\widetilde{g}_{i,\varepsilon}:=\begin{cases}
\begin{array}{cc}
g_{i} & \mathrm{if}\,\,\delta_{i}=0,\\
g_{i,\varepsilon} & \mathrm{if}\,\,\delta_{i}>0,
\end{array} & \forall\,\,i\in\{1,\ldots,n\}\end{cases}.\label{eq:g_eps_tilde}
\end{equation}

\noindent \end{brem}
\begin{lem}
\noindent \label{lem:ellipticity}Let $p\geq2$, $\varepsilon\in(0,\min\,\{1,\inf\Delta^{+}\})$
and $\xi\in\mathbb{R}^{n}$. Then, for every $\zeta\in\mathbb{R}^{n}$
we have
\begin{equation}
\varepsilon\,(1+\vert\xi\vert^{2})^{\frac{p-2}{2}}\,\vert\zeta\vert^{2}\,\leq\,\langle D^{2}F_{\varepsilon}(\xi)\,\zeta,\zeta\rangle\,\leq\,(1+\varepsilon)\,(p-1)\,(1+\vert\xi\vert^{2})^{\frac{p-2}{2}}\,\vert\zeta\vert^{2}.\label{eq:ell-bounds}
\end{equation}
\end{lem}

\noindent \begin{proof}[\bfseries{Proof}]For $p>2$, a straightforward
computation reveals that
\[
D^{2}F_{\varepsilon}(\xi)\,=\,\mathrm{diag}\left(g_{1}''(\xi_{1}),\ldots,g_{n}''(\xi_{n})\right)\,+\,\varepsilon\,D^{2}G(\xi)\,,
\]
where, for every $i\in\{1,\ldots,n\}$, 
\[
g_{i}''(s)=(p-1)\,(\vert s\vert-\delta_{i})_{+}^{p-2},\,\,\,\,\,\,\,\,\,\,s\in\mathbb{R}\,,
\]
and
\[
D^{2}G(\xi)=(1+\vert\xi\vert^{2})^{\frac{p-4}{2}}\,[(1+\vert\xi\vert^{2})\,\mathbb{I}+(p-2)\,\xi\otimes\xi]\,.
\]
Thus, for every $\zeta\in\mathbb{R}^{n}$ we get\begin{align}\label{eq:ellipticity}
\langle D^{2}F_{\varepsilon}(\xi)\,\zeta,\zeta\rangle&=\sum_{i=1}^{n}g_{i}''(\xi_{i})\,\zeta_{i}^{2}\,+\,\varepsilon\,(1+\vert\xi\vert^{2})^{\frac{p-2}{2}}\,\vert\zeta\vert^{2}\,+\,\varepsilon\,(p-2)\,(1+\vert\xi\vert^{2})^{\frac{p-4}{2}}\sum_{i,j=1}^{n}\xi_{i}\,\xi_{j}\,\zeta_{i}\,\zeta_{j}\nonumber\\
&=(p-1)\sum_{i=1}^{n}(\vert\xi_{i}\vert-\delta_{i})_{+}^{p-2}\,\zeta_{i}^{2}\,+\,\varepsilon\,(1+\vert\xi\vert^{2})^{\frac{p-2}{2}}\,\vert\zeta\vert^{2}\,+\,\varepsilon\,(p-2)\,(1+\vert\xi\vert^{2})^{\frac{p-4}{2}}\,\langle\xi,\zeta\rangle^{2}.
\end{align}Using the Cauchy-Schwarz inequality, from \eqref{eq:ellipticity}
we obtain\begin{align}\label{eq:ellipticity_02}
\langle D^{2}F_{\varepsilon}(\xi)\,\zeta,\zeta\rangle&\leq(p-1)\,(1+\vert\xi\vert^{2})^{\frac{p-2}{2}}\,\vert\zeta\vert^{2}\,+\,\varepsilon\,(1+\vert\xi\vert^{2})^{\frac{p-2}{2}}\,\vert\zeta\vert^{2}\,+\,\varepsilon\,(p-2)\,(1+\vert\xi\vert^{2})^{\frac{p-4}{2}}\,\vert\xi\vert^{2}\vert\zeta\vert^{2}\nonumber\\
&\leq(1+\varepsilon)\,(p-1)\,(1+\vert\xi\vert^{2})^{\frac{p-2}{2}}\,\vert\zeta\vert^{2}.
\end{align}For the derivation of the lower bound, it is sufficient to observe
that the first and third terms in the right-hand side of \eqref{eq:ellipticity}
are non-negative.\\
$\hspace*{1em}$In the case $p=2$, we replace each $g_{i}$ with
the function $\widetilde{g}_{i,\varepsilon}\in C^{2}(\mathbb{R})$
defined by (\ref{eq:g_eps}) and (\ref{eq:g_eps_tilde}). Noting that
$0\leq\widetilde{g}_{i,\varepsilon}''\leq1$ for any $i\in\{1,\ldots,n\}$,
when $p=2$ we immediately have 
\begin{equation}
\varepsilon\,\vert\zeta\vert^{2}\,\leq\,\langle D^{2}F_{\varepsilon}(\xi)\,\zeta,\zeta\rangle\,=\sum_{i=1}^{n}\widetilde{g}_{i,\varepsilon}''(\xi_{i})\,\zeta_{i}^{2}\,+\,\varepsilon\,\vert\zeta\vert^{2}\,\leq(1+\varepsilon)\,\vert\zeta\vert^{2}.\label{eq:ell_F_quad}
\end{equation}
This completes the proof.\end{proof}

\noindent $\hspace*{1em}$Now, for every $\varepsilon\in(0,\min\,\{1,\inf\Delta^{+}\})$,
we consider a local weak solution $u_{\varepsilon}\in L_{loc}^{p}(I;W_{loc}^{1,p}(\Omega))$
of the equation 
\[
\partial_{t}v\,=\,\mathrm{div}\,[D_{\xi}F_{\varepsilon}(Dv)]\,\,\,\,\,\,\,\,\,\,\mathrm{in}\,\,\,\Omega\times I.
\]
This means that $u_{\varepsilon}$ verifies
\begin{equation}
\iint_{\Omega\times I}\left(u_{\varepsilon}\,\partial_{t}\varphi\,-\langle D_{\xi}F_{\varepsilon}(Du_{\varepsilon}),D\varphi\rangle\right)dx\,dt\,=\,0\,,\,\,\,\,\,\,\,\,\mathrm{for\,\,every\,\,}\varphi\in C_{0}^{\infty}(\Omega\times I)\,.\label{eq:approx}
\end{equation}
Since $F_{\varepsilon}$ belongs to $C^{2}(\mathbb{R}^{n})$ and satisfies
(\ref{eq:ell-bounds}), we can rely on the classical regularity theory
for quasilinear parabolic equations, see e.g. \cite[Theorem 5.1, Chapter VIII]{DiBe}
and \cite[Lemma 3.1]{BoDuMa}, to get: 
\[
Du_{\varepsilon}\in L_{loc}^{\infty}(\Omega\times I,\mathbb{R}^{n})\,\,\,\,\,\,\,\,\,\,\,\,\mathrm{and}\,\,\,\,\,\,\,\,\,\,\,\,u_{\varepsilon}\in L_{loc}^{2}(I;W_{loc}^{2,2}(\Omega))\,.
\]
For convenience of notation, from now on we drop the index $\varepsilon\in(0,\min\,\{1,\inf\Delta^{+}\})$
and simply write $u$ and $F$ in place of $u_{\varepsilon}$ and
$F_{\varepsilon}$, unless otherwise specified.

\subsection{Caccioppoli-type inequalities\label{subsec:Caccioppoli}}

$\hspace*{1em}$The first technical tools in the proof of Theorem
\ref{thm:main} are the following standard and nonstandard Caccioppoli
inequalities, established in \cite[Lemmas 3.1 and 3.2]{BBLV}.
\begin{lem}[\textbf{Standard Caccioppoli inequality}]
\label{lem:standard}Let $\eta\in C_{0}^{\infty}(\Omega)$ and $\chi\in C_{0}^{\infty}((T_{0},T_{1}])$
be two non-negative functions, with $\chi$ non-decreasing. Let $h:\mathbb{R}\to\mathbb{R}$
be a $C^{1}$ convex non-negative function. Then, for almost every
$\tau\in I$ and every $j\in\{1,\ldots,n\}$, we have\begin{align*}
&\chi(\tau)\int_{\Omega\times\{\tau\}}h^{2}(u_{x_{j}})\,\eta^{2}\,dx+\iint_{\Omega\times(T_{0},\tau)}\langle D^{2}F(Du)\,Dh(u_{x_{j}}),Dh(u_{x_{j}})\rangle\,\chi\,\eta^{2}\,dx\,dt\nonumber\\
&\,\,\,\,\,\,\,\leq\iint_{\Omega\times(T_{0},\tau)}(\partial_{t}\chi)\,\eta^{2}\,h^{2}(u_{x_{j}})\,dx\,dt\,+\,4\iint_{\Omega\times(T_{0},\tau)}\langle D^{2}F(Du)\,D\eta,D\eta\rangle\,h^{2}(u_{x_{j}})\,\chi\,dx\,dt\,.
\end{align*}
\end{lem}

\noindent $\hspace*{1em}$While the standard Caccioppoli inequality
in Lemma \ref{lem:standard} provides an estimate for suitable convex
functions of a single spatial derivative $u_{x_{j}}$, it is not sufficient
in the orthotropic setting due to the strong degeneracy of the operator.
To overcome this issue, we need a more refined, nonstandard estimate
involving two different components of $Du$. This estimate is the
content of the next lemma, which is referred to as a \textit{weird}
Caccioppoli inequality in \cite{BBLV}. The key idea is to mix together
the components of the spatial gradient $Du$ with respect to two orthogonal
directions, in such a way as to compensate for the severe degeneracy
of the orthotropic operator, as already explained in Subsection \ref{subsec:Strategy}
and in \cite[Section 1.3]{BBLV}.
\begin{lem}[\textbf{Weird Caccioppoli inequality}]
\label{lem:weird1}Let $\eta\in C_{0}^{\infty}(\Omega)$ and $\chi\in C_{0}^{\infty}((T_{0},T_{1}])$
be two non-negative functions, with $\chi$ non-decreasing. Let $\Phi:[0,\infty)\to\mathbb{R}$
and $\Psi:[0,\infty)\to\mathbb{R}$ be two $C^{1}$ non-decreasing
and non-negative convex functions. Then, for almost every $\tau\in I$,
every $j,k\in\{1,\ldots,n\}$ and every $\alpha\in[0,1]$, we have\begin{align*}
&\chi(\tau)\int_{\Omega\times\{\tau\}}\Phi(u_{x_{j}}^{2})\,\Psi(u_{x_{k}}^{2})\,\eta^{2}\,dx\,+\iint_{\Omega\times(T_{0},\tau)}\langle D^{2}F(Du)\,Du_{x_{j}},Du_{x_{j}}\rangle\,\Phi'(u_{x_{j}}^{2})\,\Psi(u_{x_{k}}^{2})\,\chi\,\eta^{2}\,dx\,dt\nonumber\\
&\,\,\,\,\,\,\leq\iint_{\Omega\times(T_{0},\tau)}\Phi(u_{x_{j}}^{2})\,\Psi(u_{x_{k}}^{2})\,(\partial_{t}\chi)\,\eta^{2}\,dx\,dt\nonumber\\
&\,\,\,\,\,\,\,\,\,\,\,\,+\,4\iint_{\Omega\times(T_{0},\tau)}\langle D^{2}F(Du)\,D\eta,D\eta\rangle\big[u_{x_{j}}^{2}\,\Phi'(u_{x_{j}}^{2})\,\Psi(u_{x_{k}}^{2})+u_{x_{k}}^{2}\,\Phi(u_{x_{j}}^{2})\,\Psi'(u_{x_{k}}^{2})\big]\chi\,dx\,dt\nonumber\\
&\,\,\,\,\,\,\,\,\,\,\,\,+\,8\left(\iint_{\Omega\times(T_{0},\tau)}\langle D^{2}F(Du)\,Du_{x_{j}},Du_{x_{j}}\rangle\,u_{x_{j}}^{2}\,[\Phi'(u_{x_{j}}^{2})]^{2}\,[\Psi'(u_{x_{k}}^{2})]^{\alpha}\,\chi\,\eta^{2}\,dx\,dt\right)^{\frac{1}{2}}\nonumber\\
&\,\,\,\,\,\,\,\,\,\,\,\,\times\left(\iint_{\Omega\times(T_{0},\tau)}\left[\frac{1}{4}\,(\partial_{t}\chi)\,\eta^{2}+\langle D^{2}F(Du)\,D\eta,D\eta\rangle\,\chi\right]\vert u_{x_{k}}\vert^{2\alpha}\,[\Psi(u_{x_{k}}^{2})]^{2-\alpha}\,dx\,dt\right)^{\frac{1}{2}}.
\end{align*}
\end{lem}

\section{A quantitative $L^{\infty}$ bound for $Du_{\varepsilon}$\label{sec:Moser} }

\noindent $\hspace*{1em}$In this section, we establish a uniform
$L^{\infty}$ estimate for $Du_{\varepsilon}$, working with the anisotropic
parabolic cylinders defined in (\ref{eq:anis_cylin}). As in the previous
section, we will drop the index $\varepsilon\in(0,\min\,\{1,\inf\Delta^{+}\})$
and simply write $u$ and $F$ in place of $u_{\varepsilon}$ and
$F_{\varepsilon}$, respectively. Moreover, we set 
\begin{equation}
\delta:=\,1+\max\,\{\delta_{i}:i=1,\ldots,n\}\,.\label{eq:delta}
\end{equation}

\begin{prop}
\noindent \label{prop:quantitative}Let $n\geq2$ and $p\geq2$. Moreover,
let $\delta$ be defined according to $(\ref{eq:delta})$. Then, there
exist constants $\vartheta=\vartheta(n)>2$ and $C=C(n,p,\delta)>1$
such that, for every $\varepsilon\in(0,\min\,\{1,\inf\Delta^{+}\})$
and for every $Q_{r}(x_{0},t_{0})\subset Q_{R}(x_{0},t_{0})\Subset\Omega\times I$
with $R\leq1$, we have 
\begin{equation}
\Vert Du_{\varepsilon}\Vert_{L^{\infty}(Q_{r}(x_{0},t_{0}))}\,\leq\,\frac{C}{(R-r)^{\vartheta p}}\left[1+\left(\iint_{Q_{R}(x_{0},t_{0})}\vert Du_{\varepsilon}\vert^{p}\,dx\,dt\right)^{\frac{1}{2}}\right].\label{eq:uniform}
\end{equation}
\end{prop}

\noindent \begin{proof}[\bfseries{Proof}]For simplicity, we limit
ourselves to the case $n\geq3$. This allows us to apply the Sobolev
inequality valid for every $v\in W_{0}^{1,2}(\Omega)$ 
\[
\Vert v\Vert_{L^{2^{*}}(\Omega)}\,\leq\,C_{n}\,\Vert Dv\Vert_{L^{2}(\Omega)}\,\,\,\,\,\,\,\,\,\,\,\,\mathrm{with}\,\,2^{*}=\frac{2n}{n-2}\,.
\]

\noindent Here $C_{n}$ denotes a constant depending only on $n$.
The case $n=2$ requires only minor modifications, whose details are
left to the reader.

\noindent $\hspace*{1em}$As the proof is rather intricate, we divide
it into several steps for clarity. The first two steps follow the
proof of the analogous Proposition 4.1 in \cite{BBLV}, where, however,
$\delta_{i}=0$ for every $i\in\{1,\ldots,n\}$. For the sake of completeness,
we include them here as well. \medskip{}

\noindent \textbf{Step 1: the choices of $\Phi$ and $\Psi$.} We
apply Lemma \ref{lem:weird1} with the following choices:
\[
\Phi(t)=t^{s}\,\,\,\,\,\,\,\,\,\,\,\,\mathrm{and}\,\,\,\,\,\,\,\,\,\,\,\,\Psi(t)=t^{m},\,\,\,\,\,\,\,\,\,\,\mathrm{for}\,\,t\geq0\,,
\]
with $1\leq s\leq m$. We also take 
\[
\alpha=\begin{cases}
\begin{array}{cc}
{\displaystyle \frac{m-s}{m-1}}\in[0,1] & \mathrm{if}\,\,m>1,\vspace{2mm}\\
1\,\,\,\,\,\,\,\,\,\,\,\,\,\,\,\,\,\,\,\,\,\,\,\,\,\,\,\,\,\,\,\,\, & \mathrm{if}\,\,m=1.
\end{array}\end{cases}
\]
This yields\begin{align*}
&\chi(\tau)\int_{\Omega\times\{\tau\}}\vert u_{x_{j}}\vert^{2s}\,\vert u_{x_{k}}\vert^{2m}\,\eta^{2}\,dx\,+\,s\iint_{\Omega\times(T_{0},\tau)}\langle D^{2}F(Du)\,Du_{x_{j}},Du_{x_{j}}\rangle\,\vert u_{x_{j}}\vert^{2s-2}\,\vert u_{x_{k}}\vert^{2m}\,\chi\,\eta^{2}\,dx\,dt\\
&\,\,\,\,\,\,\,\leq\iint_{\Omega\times(T_{0},\tau)}\vert u_{x_{j}}\vert^{2s}\,\vert u_{x_{k}}\vert^{2m}\,(\partial_{t}\chi)\,\eta^{2}\,dx\,dt\\
&\,\,\,\,\,\,\,\,\,\,\,\,\,\,+\,4\iint_{\Omega\times(T_{0},\tau)}\langle D^{2}F(Du)\,D\eta,D\eta\rangle\big[s\,\vert u_{x_{j}}\vert^{2s}\,\vert u_{x_{k}}\vert^{2m}+m\,\vert u_{x_{k}}\vert^{2m}\,\vert u_{x_{j}}\vert^{2s}\big]\,\chi\,dx\,dt\\
&\,\,\,\,\,\,\,\,\,\,\,\,\,\,+\,8\,s\,m^{\frac{\alpha}{2}}\left(\iint_{\Omega\times(T_{0},\tau)}\langle D^{2}F(Du)\,Du_{x_{j}},Du_{x_{j}}\rangle\,\vert u_{x_{j}}\vert^{4s-2}\,\vert u_{x_{k}}\vert^{2m-2s}\,\chi\,\eta^{2}\,dx\,dt\right)^{\frac{1}{2}}\\
&\,\,\,\,\,\,\,\,\,\,\,\,\,\,\times\left(\iint_{\Omega\times(T_{0},\tau)}\left[\frac{1}{4}\,(\partial_{t}\chi)\,\eta^{2}+\langle D^{2}F(Du)\,D\eta,D\eta\rangle\,\chi\right]\vert u_{x_{k}}\vert^{2(s+m)}\,dx\,dt\right)^{\frac{1}{2}}.
\end{align*}On the product of the last two integrals, we apply Young's inequality
in the form 
\[
ab\leq a^{2}+\frac{b^{2}}{4}\,.
\]
Thus we obtain\begin{align*}
&\chi(\tau)\int_{\Omega\times\{\tau\}}\vert u_{x_{j}}\vert^{2s}\,\vert u_{x_{k}}\vert^{2m}\,\eta^{2}\,dx\,+\iint_{\Omega\times(T_{0},\tau)}\langle D^{2}F(Du)\,Du_{x_{j}},Du_{x_{j}}\rangle\,\vert u_{x_{j}}\vert^{2s-2}\,\vert u_{x_{k}}\vert^{2m}\,\chi\,\eta^{2}\,dx\,dt\nonumber\\
&\,\,\,\,\,\,\,\leq\iint_{\Omega\times(T_{0},\tau)}\vert u_{x_{j}}\vert^{2s}\,\vert u_{x_{k}}\vert^{2m}\,(\partial_{t}\chi)\,\eta^{2}\,dx\,dt\nonumber\\
&\,\,\,\,\,\,\,\,\,\,\,\,\,\,+\,4\,(s+m)\iint_{\Omega\times(T_{0},\tau)}\langle D^{2}F(Du)\,D\eta,D\eta\rangle\,\vert u_{x_{j}}\vert^{2s}\,\vert u_{x_{k}}\vert^{2m}\,\chi\,dx\,dt\nonumber\\
&\,\,\,\,\,\,\,\,\,\,\,\,\,\,+\,16\,s^{2}\,m\iint_{\Omega\times(T_{0},\tau)}\left[\frac{1}{4}\,(\partial_{t}\chi)\,\eta^{2}+\langle D^{2}F(Du)\,D\eta,D\eta\rangle\,\chi\right]\vert u_{x_{k}}\vert^{2(s+m)}\,dx\,dt\nonumber\\
&\,\,\,\,\,\,\,\,\,\,\,\,\,\,+\iint_{\Omega\times(T_{0},\tau)}\langle D^{2}F(Du)\,Du_{x_{j}},Du_{x_{j}}\rangle\,\vert u_{x_{j}}\vert^{4s-2}\,\vert u_{x_{k}}\vert^{2m-2s}\,\chi\,\eta^{2}\,dx\,dt\,,
\end{align*}where we have also used that $s\geq1$ in the left-hand side and $m^{\alpha}\leq m$
in the right-hand side. By Young's inequality again, we can estimate
\[
\vert u_{x_{j}}\vert^{2s}\,\vert u_{x_{k}}\vert^{2m}\,\leq\,\vert u_{x_{j}}\vert^{2(s+m)}+\vert u_{x_{k}}\vert^{2(s+m)}\,.
\]
This finally gives\begin{align}\label{eq:step1_stima02}
&\chi(\tau)\int_{\Omega\times\{\tau\}}\vert u_{x_{j}}\vert^{2s}\,\vert u_{x_{k}}\vert^{2m}\,\eta^{2}\,dx\,+\iint_{\Omega\times(T_{0},\tau)}\langle D^{2}F(Du)\,Du_{x_{j}},Du_{x_{j}}\rangle\,\vert u_{x_{j}}\vert^{2s-2}\,\vert u_{x_{k}}\vert^{2m}\,\chi\,\eta^{2}\,dx\,dt\nonumber\\
&\,\leq16\,(s+m+s^{2}\,m)\iint_{\Omega\times(T_{0},\tau)}\left[(\partial_{t}\chi)\,\eta^{2}+\langle D^{2}F(Du)\,D\eta,D\eta\rangle\,\chi\right]\left[\vert u_{x_{j}}\vert^{2(s+m)}+\vert u_{x_{k}}\vert^{2(s+m)}\right]dx\,dt\nonumber\\
&\,\,\,\,\,\,\,\,+\iint_{\Omega\times(T_{0},\tau)}\langle D^{2}F(Du)\,Du_{x_{j}},Du_{x_{j}}\rangle\,\vert u_{x_{j}}\vert^{4s-2}\,\vert u_{x_{k}}\vert^{2m-2s}\,\chi\,\eta^{2}\,dx\,dt\,.
\end{align}\\
\textbf{Step 2: the staircase.} Let $\ell_{0}\in\mathbb{N\,}\backslash\,\{0\}$
and set $M=2^{\ell_{0}}-1$. We define the two families of indices
\[
s_{\ell}\,=\,2^{\ell}\,\,\,\,\,\,\,\,\,\,\,\,\mathrm{and}\,\,\,\,\,\,\,\,\,\,\,\,m_{\ell}\,=\,M+1-2^{\ell},\,\,\,\,\,\,\,\,\,\,\mathrm{for}\,\,\ell\in\{0,\ldots,\ell_{0}\}\,.
\]
By construction, for every $\ell\in\{0,\ldots,\ell_{0}-1\}$ we have
\[
s_{\ell}+m_{\ell}\,=\,M+1,\,\,\,\,\,\,\,\,\,\,\,\,2s_{\ell}-1\,=\,s_{\ell+1}-1\,\,\,\,\,\,\,\,\,\,\,\,\mathrm{and}\,\,\,\,\,\,\,\,\,\,\,\,m_{\ell}-s_{\ell}\,=\,m_{\ell+1}\,.
\]
We also use that $s_{\ell}+m_{\ell}+s_{\ell}^{2}\,m_{\ell}\leq2(M+1)^{3}$.
Then, inequality \eqref{eq:step1_stima02} written for $s=s_{\ell}$
and $m=m_{\ell}$, with $0\leq\ell\leq\ell_{0}-1$, gives\begin{align*}
&\chi(\tau)\int_{\Omega\times\{\tau\}}\vert u_{x_{j}}\vert^{2s_{\ell}}\,\vert u_{x_{k}}\vert^{2m_{\ell}}\,\eta^{2}\,dx\,+\iint_{\Omega\times(T_{0},\tau)}\langle D^{2}F(Du)\,Du_{x_{j}},Du_{x_{j}}\rangle\,\vert u_{x_{j}}\vert^{2s_{\ell}-2}\,\vert u_{x_{k}}\vert^{2m_{\ell}}\,\chi\,\eta^{2}\,dx\,dt\nonumber\\
&\,\,\,\,\,\leq\,32\,(M+1)^{3}\iint_{\Omega\times(T_{0},\tau)}\left[(\partial_{t}\chi)\,\eta^{2}+\langle D^{2}F(Du)\,D\eta,D\eta\rangle\,\chi\right]\left[\vert u_{x_{j}}\vert^{2(M+1)}+\vert u_{x_{k}}\vert^{2(M+1)}\right]dx\,dt\nonumber\\
&\,\,\,\,\,\,\,\,\,\,\,\,+\iint_{\Omega\times(T_{0},\tau)}\langle D^{2}F(Du)\,Du_{x_{j}},Du_{x_{j}}\rangle\,\vert u_{x_{j}}\vert^{2s_{\ell+1}-\,2}\,\vert u_{x_{k}}\vert^{2m_{\ell+1}}\,\chi\,\eta^{2}\,dx\,dt\,.
\end{align*}By summing with respect to $\ell$ from $0$ to $\ell_{0}-1$, and
erasing the common terms on both sides, we get\begin{align*}
&\chi(\tau)\,\sum_{\ell=0}^{\ell_{0}-1}\int_{\Omega\times\{\tau\}}\vert u_{x_{j}}\vert^{2s_{\ell}}\,\vert u_{x_{k}}\vert^{2m_{\ell}}\,\eta^{2}\,dx\,+\iint_{\Omega\times(T_{0},\tau)}\langle D^{2}F(Du)\,Du_{x_{j}},Du_{x_{j}}\rangle\,\vert u_{x_{k}}\vert^{2M}\,\chi\,\eta^{2}\,dx\,dt\\
&\,\,\,\,\,\,\,\leq\,260\,M^{3}\,\ell_{0}\iint_{\Omega\times(T_{0},\tau)}\left[(\partial_{t}\chi)\,\eta^{2}+\langle D^{2}F(Du)\,D\eta,D\eta\rangle\,\chi\right]\left[\vert u_{x_{j}}\vert^{2(M+1)}+\vert u_{x_{k}}\vert^{2(M+1)}\right]dx\,dt\\
&\,\,\,\,\,\,\,\,\,\,\,\,\,\,+\iint_{\Omega\times(T_{0},\tau)}\langle D^{2}F(Du)\,Du_{x_{j}},Du_{x_{j}}\rangle\,\vert u_{x_{j}}\vert^{2M}\,\chi\,\eta^{2}\,dx\,dt\,.
\end{align*}Now, in order to estimate the last term on the right-hand side of
the previous inequality, we use Lemma \ref{lem:standard} with the
choice 
\[
h(y)\,=\,\frac{\vert y\vert^{M+1}}{M+1}\,,\,\,\,\,\,\,\,\,\,\,\,\,y\in\mathbb{R}\,.
\]
We thus obtain\begin{align*}
&\iint_{\Omega\times(T_{0},\tau)}\langle D^{2}F(Du)\,Du_{x_{j}},Du_{x_{j}}\rangle\,\vert u_{x_{j}}\vert^{2M}\,\chi\,\eta^{2}\,dx\,dt\\
&\,\,\,\,\,\,\,=\iint_{\Omega\times(T_{0},\tau)}\langle D^{2}F(Du)\,Dh(u_{x_{j}}),Dh(u_{x_{j}})\rangle\,\chi\,\eta^{2}\,dx\,dt\\
&\,\,\,\,\,\,\,\leq\,\frac{4}{(M+1)^{2}}\iint_{\Omega\times(T_{0},\tau)}[(\partial_{t}\chi)\,\eta^{2}+\langle D^{2}F(Du)\,D\eta,D\eta\rangle\,\chi]\,\vert u_{x_{j}}\vert^{2(M+1)}\,dx\,dt\,.
\end{align*}Combining the two previous estimates and using the fact that 
\[
M\,=\,2^{\ell_{0}}-1\,\geq\,\ell_{0}\,\geq\,1\,\,\,\,\,\,\,\,\,\,\,\,\mathrm{for\,\,every\,\,}\ell_{0}\in\mathbb{N\,\backslash\,}\{0\},
\]
we find\begin{align}\label{eq:step3_stima01}
&\chi(\tau)\,\sum_{\ell=0}^{\ell_{0}-1}\int_{\Omega\times\{\tau\}}\vert u_{x_{j}}\vert^{2s_{\ell}}\,\vert u_{x_{k}}\vert^{2m_{\ell}}\,\eta^{2}\,dx\,+\iint_{\Omega\times(T_{0},\tau)}\langle D^{2}F(Du)\,Du_{x_{j}},Du_{x_{j}}\rangle\,\vert u_{x_{k}}\vert^{2M}\,\chi\,\eta^{2}\,dx\,dt\nonumber\\
&\,\,\,\,\,\,\,\leq\,261\,M^{4}\iint_{\Omega\times(T_{0},\tau)}\left[(\partial_{t}\chi)\,\eta^{2}+\langle D^{2}F(Du)\,D\eta,D\eta\rangle\,\chi\right]\left[\vert u_{x_{j}}\vert^{2(M+1)}+\vert u_{x_{k}}\vert^{2(M+1)}\right]dx\,dt\,.
\end{align}\\
\textbf{Step 3: weak ellipticity and boundedness of $D^{2}F$ for
$p>2$.} From now on, we shall assume that $p>2$, unless otherwise
specified.\\
$\hspace*{1em}$From \eqref{eq:ellipticity} and \eqref{eq:ellipticity_02},
it follows that
\[
(p-1)\sum_{i=1}^{n}(\vert\xi_{i}\vert-\delta_{i})_{+}^{p-2}\,\zeta_{i}^{2}\,\leq\langle D^{2}F(\xi)\,\zeta,\zeta\rangle\leq\,c(p)\,(1+\vert\xi\vert^{p-2})\,\vert\zeta\vert^{2}\,\,\,\,\,\,\,\,\,\,\,\,\mathrm{for\,\,every\,\,}\xi,\zeta\in\mathbb{R}^{n}.
\]
Inserting these estimates into \eqref{eq:step3_stima01}, one gets\begin{align*}
&\chi(\tau)\,\sum_{\ell=0}^{\ell_{0}-1}\int_{\Omega\times\{\tau\}}\vert u_{x_{j}}\vert^{2s_{\ell}}\,\vert u_{x_{k}}\vert^{2m_{\ell}}\,\eta^{2}\,dx\\
&\,\,\,\,\,\,\,+\,(p-1)\,\sum_{i=1}^{n}\iint_{\Omega\times(T_{0},\tau)}(\vert u_{x_{i}}\vert-\delta_{i})_{+}^{p-2}\,u_{x_{i}x_{j}}^{2}\,\vert u_{x_{k}}\vert^{2M}\,\chi\,\eta^{2}\,dx\,dt\\
&\,\,\,\,\,\,\,\,\,\,\,\,\,\,\leq\,c(p)\,M^{4}\iint_{\Omega\times(T_{0},\tau)}[(\partial_{t}\chi)\,\eta^{2}+\chi\,(1+\vert Du\vert^{p-2})\,\vert D\eta\vert^{2}]\,[\vert u_{x_{j}}\vert^{2(M+1)}+\vert u_{x_{k}}\vert^{2(M+1)}]\,dx\,dt\,.
\end{align*}We now consider the second term on the left-hand side. By keeping
in the sum only the term with $i=k$ and dropping the others, we obtain\begin{align*}
&\sum_{i=1}^{n}\iint_{\Omega\times(T_{0},\tau)}(\vert u_{x_{i}}\vert-\delta_{i})_{+}^{p-2}\,u_{x_{i}x_{j}}^{2}\,\vert u_{x_{k}}\vert^{2M}\,\chi\,\eta^{2}\,dx\,dt\\
&\,\,\,\,\,\,\,\geq\iint_{\Omega\times(T_{0},\tau)}(\vert u_{x_{k}}\vert-\delta_{k})_{+}^{p-2}\,u_{x_{k}x_{j}}^{2}\,\vert u_{x_{k}}\vert^{2M}\,\chi\,\eta^{2}\,dx\,dt\,.
\end{align*}Note that, by Young's inequality, one has\begin{align*}
&\iint_{\Omega\times(T_{0},\tau)}\left|[(\vert u_{x_{k}}\vert-\delta_{k})_{+}^{\frac{p}{2}}\,\vert u_{x_{k}}\vert^{M}]_{x_{j}}\right|^{2}\chi\,\eta^{2}\,dx\,dt\\
&\,\,\,\,\,\,\,\leq\,2\iint_{\Omega\times(T_{0},\tau)}\left|[(\vert u_{x_{k}}\vert-\delta_{k})_{+}^{\frac{p}{2}}]_{x_{j}}\right|^{2}\vert u_{x_{k}}\vert^{2M}\,\chi\,\eta^{2}\,dx\,dt\\
&\,\,\,\,\,\,\,\,\,\,\,\,\,\,+\,2\iint_{\Omega\times(T_{0},\tau)}(\vert u_{x_{k}}\vert-\delta_{k})_{+}^{p}\left|[\vert u_{x_{k}}\vert^{M}]_{x_{j}}\right|^{2}\chi\,\eta^{2}\,dx\,dt\\
&\,\,\,\,\,\,\,\leq\left(\frac{p^{2}}{2}\,+2M^{2}\right)\iint_{\Omega\times(T_{0},\tau)}(\vert u_{x_{k}}\vert-\delta_{k})_{+}^{p-2}\,u_{x_{k}x_{j}}^{2}\vert u_{x_{k}}\vert^{2M}\,\chi\,\eta^{2}\,dx\,dt\\
&\,\,\,\,\,\,\,\leq\,p^{2}M^{2}\iint_{\Omega\times(T_{0},\tau)}(\vert u_{x_{k}}\vert-\delta_{k})_{+}^{p-2}\,u_{x_{k}x_{j}}^{2}\vert u_{x_{k}}\vert^{2M}\,\chi\,\eta^{2}\,dx\,dt\,,
\end{align*}where, in the last line, we have used that $M\geq1$ and $\frac{p^{2}}{2}>2$.
Then, combining the three previous estimates and summing over $j\in\{1,\ldots,n\}$
the resulting inequality, we get\begin{align*}
&\chi(\tau)\,\sum_{\ell=0}^{\ell_{0}-1}\int_{\Omega\times\{\tau\}}\sum_{j=1}^{n}\,\vert u_{x_{j}}\vert^{2s_{\ell}}\,\vert u_{x_{k}}\vert^{2m_{\ell}}\,\eta^{2}\,dx\\
&\,\,\,\,\,\,\,+\,\frac{p-1}{p^{2}M^{2}}\iint_{\Omega\times(T_{0},\tau)}\left|D[(\vert u_{x_{k}}\vert-\delta_{k})_{+}^{\frac{p}{2}}\,\vert u_{x_{k}}\vert^{M}]\right|^{2}\chi\,\eta^{2}\,dx\,dt\\
&\,\,\,\,\,\,\,\,\,\,\,\,\,\,\leq\,c(p)\,M^{4}\iint_{\Omega\times(T_{0},\tau)}[(\partial_{t}\chi)\,\eta^{2}+\chi\,(1+\vert Du\vert^{p-2})\,\vert D\eta\vert^{2}]\left[\sum_{j=1}^{n}\,\vert u_{x_{j}}\vert^{2M+2}+n\,\vert u_{x_{k}}\vert^{2M+2}\right]dx\,dt\,.
\end{align*} We now use that $p-1<p^{2}M^{2}$ and add the term 
\[
\iint_{\Omega\times(T_{0},\tau)}(\vert u_{x_{k}}\vert-\delta_{k})_{+}^{p}\,\vert u_{x_{k}}\vert^{2M}\,\chi\,\vert D\eta\vert^{2}\,dx\,dt
\]
to both sides of the preceding inequality. With some algebraic manipulations,
this gives\begin{align*}
&\chi(\tau)\,\sum_{\ell=0}^{\ell_{0}-1}\int_{\Omega\times\{\tau\}}\sum_{j=1}^{n}\,\vert u_{x_{j}}\vert^{2s_{\ell}}\,\vert u_{x_{k}}\vert^{2m_{\ell}}\,\eta^{2}\,dx\,+\iint_{\Omega\times(T_{0},\tau)}\left|D[(\vert u_{x_{k}}\vert-\delta_{k})_{+}^{\frac{p}{2}}\,\vert u_{x_{k}}\vert^{M}\,\eta]\right|^{2}\chi\,dx\,dt\\
&\,\,\,\,\,\,\,\leq\,c(p)\,M^{6}\iint_{\Omega\times(T_{0},\tau)}[(\partial_{t}\chi)\,\eta^{2}+\chi\,(1+\vert Du\vert^{p-2})\,\vert D\eta\vert^{2}]\left[\sum_{j=1}^{n}\,\vert u_{x_{j}}\vert^{2M+2}+n\,\vert u_{x_{k}}\vert^{2M+2}\right]dx\,dt\\
&\,\,\,\,\,\,\,\,\,\,\,\,\,\,+\,2\iint_{\Omega\times(T_{0},\tau)}(\vert u_{x_{k}}\vert-\delta_{k})_{+}^{p}\,\vert u_{x_{k}}\vert^{2M}\,\chi\,\vert D\eta\vert^{2}\,dx\,dt\\
&\,\,\,\,\,\,\,\leq\,c(p)\,M^{6}\iint_{\Omega\times(T_{0},\tau)}[(\partial_{t}\chi)\,\eta^{2}+\chi\,(1+\vert Du\vert^{p-2})\,\vert D\eta\vert^{2}]\left[\sum_{j=1}^{n}\,\vert u_{x_{j}}\vert^{2M+2}+n\,\vert u_{x_{k}}\vert^{2M+2}\right]dx\,dt\\
&\,\,\,\,\,\,\,\,\,\,\,\,\,\,+\,c(p)\,M^{6}\iint_{\Omega\times(T_{0},\tau)}\vert u_{x_{k}}\vert^{p+2M}\,\chi\,\vert D\eta\vert^{2}\,dx\,dt\,.
\end{align*}By using the Sobolev inequality in the spatial variable for the second
term on the left-hand side, we obtain\begin{align*}
&\chi(\tau)\,\sum_{\ell=0}^{\ell_{0}-1}\int_{\Omega\times\{\tau\}}\sum_{j=1}^{n}\,\vert u_{x_{j}}\vert^{2s_{\ell}}\,\vert u_{x_{k}}\vert^{2m_{\ell}}\,\eta^{2}\,dx\,+\int_{T_{0}}^{\tau}\chi\left(\int_{\Omega}(\vert u_{x_{k}}\vert-\delta_{k})_{+}^{p\,\frac{2^{*}}{2}}\,\vert u_{x_{k}}\vert^{2^{*}M}\,\eta^{2^{*}}dx\right)^{\frac{2}{2^{*}}}dt\\
&\,\,\,\,\,\,\,\leq\,c\,M^{6}\iint_{\Omega\times(T_{0},\tau)}[(\partial_{t}\chi)\,\eta^{2}+\chi\,(1+\vert Du\vert^{p-2})\,\vert D\eta\vert^{2}]\left[\sum_{j=1}^{n}\,\vert u_{x_{j}}\vert^{2M+2}+n\,\vert u_{x_{k}}\vert^{2M+2}\right]dx\,dt\\
&\,\,\,\,\,\,\,\,\,\,\,\,\,\,+\,c\,M^{6}\iint_{\Omega\times(T_{0},\tau)}\vert u_{x_{k}}\vert^{p+2M}\,\chi\,\vert D\eta\vert^{2}\,dx\,dt\,,
\end{align*}where $c$ is now a positive constant depending only on $n$ and $p$.
Finally, we sum over $k\in\{1,\ldots,n\}$ and apply Minkowski's inequality
to the second term on the left-hand side. This yields\begin{align}\label{eq:step3_stima002}
&\chi(\tau)\,\sum_{\ell=0}^{\ell_{0}-1}\int_{\Omega\times\{\tau\}}\sum_{j=1}^{n}\,\vert u_{x_{j}}\vert^{2s_{\ell}}\,\sum_{k=1}^{n}\,\vert u_{x_{k}}\vert^{2m_{\ell}}\,\eta^{2}\,dx\nonumber\\
&\,\,\,\,+\int_{T_{0}}^{\tau}\chi\left(\int_{\Omega}\left|\sum_{k=1}^{n}\,(\vert u_{x_{k}}\vert-\delta_{k})_{+}^{p}\,\vert u_{x_{k}}\vert^{2M}\right|^{\frac{2^{*}}{2}}\eta^{2^{*}}dx\right)^{\frac{2}{2^{*}}}dt\nonumber\\
&\,\,\,\,\,\,\,\,\,\,\,\leq\,c\,M^{6}\iint_{\Omega\times(T_{0},\tau)}[(\partial_{t}\chi)\,\eta^{2}+\chi\,(1+\vert Du\vert^{p-2})\,\vert D\eta\vert^{2}]\left[\sum_{j=1}^{n}\,\vert u_{x_{j}}\vert^{2M+2}+\sum_{k=1}^{n}\,\vert u_{x_{k}}\vert^{2M+2}\right]dx\,dt\nonumber\\
&\,\,\,\,\,\,\,\,\,\,\,\,\,\,\,\,\,\,+\,c\,M^{6}\iint_{\Omega\times(T_{0},\tau)}\sum_{k=1}^{n}\,\vert u_{x_{k}}\vert^{p+2M}\,\chi\,\vert D\eta\vert^{2}\,dx\,dt\,.
\end{align}We now introduce the auxiliary function 
\[
\mathcal{U}(x,t):=\,\frac{1}{2\delta}\,\,\underset{1\,\le\,i\,\leq\,n}{\max}\vert u_{x_{i}}(x,t)\vert\,,
\]
where the parameter $\delta$ is defined in (\ref{eq:delta}). A few
elementary calculations show that
\[
\sum_{k=1}^{n}(\vert u_{x_{k}}\vert-\delta_{k})_{+}^{p}\,\vert u_{x_{k}}\vert^{2M}\,\ge\,(2\delta)^{p+2M}\left(\mathcal{U}-\frac{1}{2}\right)_{+}^{p}\,\mathcal{U}^{2M}
\]
and
\[
(2\delta\mathcal{U})^{q}\,\le\,\sum_{k=1}^{n}\,\vert u_{x_{k}}\vert^{q}\,\leq\,n\,(2\delta\mathcal{U})^{q}\,\,\,\,\,\,\,\,\,\,\,\,\mathrm{for\,\,every}\,\,q\geq0\,.
\]
In particular, for $q=2$ we obtain 
\begin{equation}
2\delta\mathcal{U}\,\leq\,\vert Du\vert\,\leq\,\sqrt{n}\,2\delta\mathcal{U}\,.\label{eq:controlDu}
\end{equation}
Inserting these estimates into \eqref{eq:step3_stima002} yields \begin{align}\label{eq:step3_stima003}
&\chi(\tau)\,\ell_{0}\,(2\delta)^{2M+2}\int_{\Omega\times\{\tau\}}\mathcal{U}^{2M+2}\,\eta^{2}\,dx\,+\,(2\delta)^{p+2M}\int_{T_{0}}^{\tau}\chi\left(\int_{\Omega}\left(\mathcal{U}-\frac{1}{2}\right)_{+}^{p\,\frac{2^{*}}{2}}\,\mathcal{U}^{2^{*}M}\,\eta^{2^{*}}dx\right)^{\frac{2}{2^{*}}}dt\nonumber\\
&\,\,\,\,\,\,\,\leq\,c\,M^{6}\,(2\delta)^{p+2M}\iint_{\Omega\times(T_{0},\tau)}[(\partial_{t}\chi)\,\eta^{2}+\chi\,(1+\mathcal{U}^{p-2})\,\vert D\eta\vert^{2}]\,\mathcal{U}^{2M+2}\,dx\,dt\nonumber\\
&\,\,\,\,\,\,\,\,\,\,\,\,\,\,+\,c\,M^{6}\,(2\delta)^{p+2M}\iint_{\Omega\times(T_{0},\tau)}\mathcal{U}^{p+2M}\,\chi\,\vert D\eta\vert^{2}\,dx\,dt\nonumber\\
&\,\,\,\,\,\,\,\leq\,c\,M^{6}\,(2\delta)^{p+2M}\iint_{\Omega\times(T_{0},\tau)}[(\partial_{t}\chi)\,\eta^{2}+\chi\,\vert D\eta\vert^{2}]\,(1+\mathcal{U}^{p+2M})\,dx\,dt\,,
\end{align}where, in the first line, we have also used that $s_{\ell}+m_{\ell}=M+1$
for every $\ell\in\{0,\ldots,\ell_{0}-1\}$. Dividing both sides of
\eqref{eq:step3_stima003} by $(2\delta)^{2M+2}$ and using that $\ell_{0}\geq1$
together with $(2\delta)^{p-2}\geq2^{p-2}>1$, we get\begin{align}\label{eq:step3_stima04}
&\chi(\tau)\int_{\Omega\times\{\tau\}}\mathcal{U}^{2M+2}\,\eta^{2}\,dx\,+\int_{T_{0}}^{\tau}\chi\left(\int_{\Omega}\left(\mathcal{U}-\frac{1}{2}\right)_{+}^{p\,\frac{2^{*}}{2}}\,\mathcal{U}^{2^{*}M}\,\eta^{2^{*}}dx\right)^{\frac{2}{2^{*}}}dt\nonumber\\
&\,\,\,\,\,\,\,\leq\,CM^{6}\iint_{\Omega\times(T_{0},\tau)}[(\partial_{t}\chi)\,\eta^{2}+\chi\,\vert D\eta\vert^{2}]\,(1+\mathcal{U}^{p+2M})\,dx\,dt\,,
\end{align} where $C$ is a positive constant depending only on $n$, $p$ and
$\delta$.\\
\\
\textbf{Step 4: choice of the cut-off functions.} Let $(x_{0},t_{0})\in\Omega\times I$
and $0<r<R\leq1$ be such that the cube $Q_{R}(x_{0}):=x_{0}+(-R,R)^{n}$
is compactly contained in $\Omega$. In addition, we require that
\[
(t_{0}-R^{p},t_{0})\Subset I,
\]
so that we must have $T_{0}<t_{0}<T_{1}$ and $R^{p}<t_{0}-T_{0}$.
Let $\chi:[T_{0},T_{1}]\to\mathbb{R}$ be a non-decreasing Lipschitz
continuous function such that 
\[
\chi\equiv0\,\,\,\,\mathrm{on}\,\,[T_{0},t_{0}-R^{p}],\,\,\,\,\,\,\,\,\,\,\,\,\chi\equiv1\,\,\,\,\mathrm{on}\,\,[t_{0}-r^{p},t_{0}]\,\,\,\,\,\,\,\,\,\,\,\,\mathrm{and}\,\,\,\,\,\,\,\,\,\,\,\,\partial_{t}\chi\,\leq\,\frac{\tilde{c}}{(R-r)^{p}}\,.
\]

\noindent Let $\eta\in C_{0}^{\infty}(Q_{R}(x_{0}))$ be such that
\[
0\leq\eta\leq1,\,\,\,\,\,\,\,\,\,\,\,\,\eta\equiv1\,\,\,\,\mathrm{on}\,\,Q_{r}(x_{0})\,\,\,\,\,\,\,\,\,\,\,\,\mathrm{and}\,\,\,\,\,\,\,\,\,\,\,\,\vert D\eta\vert\,\leq\,\frac{\tilde{c}}{R-r}\,.
\]
We recall the notation for the anisotropic parabolic cylinder
\[
Q_{\rho}(x_{0},t_{0}):=\,Q_{\rho}(x_{0})\times(t_{0}-\rho^{p},t_{0}),\,\,\,\,\,\,\,\,\,\,\,\,\rho>0.
\]
With this choice of $\chi$ and $\eta$, we apply estimate \eqref{eq:step3_stima04}
twice: firstly, by discarding the second term on the left-hand side
and taking the supremum over $\tau$ in the interval $(t_{0}-r^{p},t_{0})$;
secondly, by dropping the first term on the left-hand side and taking
$\tau=t_{0}$. Summing the two resulting inequalities yields\begin{align}\label{eq:step4_stima01}
&\underset{\tau\,\in\,(t_{0}\,-\,r^{p},\,t_{0})}{\sup}\int_{Q_{r}(x_{0})\times\{\tau\}}\mathcal{U}^{2M+2}\,dx\,+\int_{t_{0}\,-\,r^{p}}^{t_{0}}\left(\int_{Q_{r}(x_{0})}\left(\mathcal{U}-\frac{1}{2}\right)_{+}^{p\,\frac{2^{*}}{2}}\,\mathcal{U}^{2^{*}M}\,dx\right)^{\frac{2}{2^{*}}}dt\nonumber\\
&\,\,\,\,\,\,\,\leq\,C\,\frac{M^{6}}{(R-r)^{p}}\iint_{Q_{R}(x_{0},t_{0})}(1+\mathcal{U}^{p+2M})\,dx\,dt\,,
\end{align}where we have also used that $(R-r)^{p}\leq(R-r)^{2}$, since $R\leq1$
and $p>2$. By Hölder's inequality, we have\begin{align*}
&\iint_{Q_{r}(x_{0},t_{0})}\left(\mathcal{U}-\frac{1}{2}\right)_{+}^{p}\,\mathcal{U}^{2M\,+\,\frac{4(M+1)}{n}}\,dx\,dt\nonumber\\
&\,\,\,\,\,\,\,\leq\int_{t_{0}\,-\,r^{p}}^{t_{0}}\left(\int_{Q_{r}(x_{0})}\left(\mathcal{U}-\frac{1}{2}\right)_{+}^{p\,\frac{2^{*}}{2}}\,\mathcal{U}^{2^{*}M}\,dx\right)^{\frac{2}{2^{*}}}\left(\int_{Q_{r}(x_{0})}\mathcal{U}^{2M+2}\,dx\right)^{\frac{2}{n}}dt\nonumber\\
&\,\,\,\,\,\,\,\leq\left(\underset{\tau\,\in\,(t_{0}\,-\,r^{p},\,t_{0})}{\sup}\int_{Q_{r}(x_{0})\times\{\tau\}}\mathcal{U}^{2M+2}\,dx\right)^{\frac{2}{n}}\int_{t_{0}\,-\,r^{p}}^{t_{0}}\left(\int_{Q_{r}(x_{0})}\left(\mathcal{U}-\frac{1}{2}\right)_{+}^{p\,\frac{2^{*}}{2}}\,\mathcal{U}^{2^{*}M}\,dx\right)^{\frac{2}{2^{*}}}dt\,.
\end{align*}Combining the previous estimate with \eqref{eq:step4_stima01}, we
obtain
\begin{equation}
\iint_{Q_{r}(x_{0},t_{0})}\left(\mathcal{U}-\frac{1}{2}\right)_{+}^{p}\,\mathcal{U}^{2M\,+\,\frac{4(M+1)}{n}}\,dx\,dt\,\leq\left[C\,\frac{M^{6}}{(R-r)^{p}}\iint_{Q_{R}(x_{0},t_{0})}(1+\mathcal{U}^{p+2M})\,dx\,dt\right]^{\frac{2}{n}+1}.\label{eq:st4_stima2}
\end{equation}
We now estimate\begin{align*}
\iint_{Q_{R}(x_{0},t_{0})}(1+\mathcal{U}^{p+2M})\,dx\,dt\,&\leq\, 2\,\vert Q_{R}(x_{0},t_{0})\vert+\iint_{Q_{R}(x_{0},t_{0})\cap\{\mathcal{U}\,\geq\,1\}}\mathcal{U}^{p+2M}\,dx\,dt\\
&\leq\, 2^{n+1}+\iint_{Q_{R}(x_{0},t_{0})\cap\{\mathcal{U}\,\geq\,1\}}\mathcal{U}^{p+2M}\,dx\,dt\,,
\end{align*} where, in the last line, we have used that $R\leq1$. Observe that
on the set $\{\mathcal{U}\geq1\}$, we have $\mathcal{U}\leq2\left(\mathcal{U}-\frac{1}{2}\right)_{+}$.
Hence,\begin{align}\label{eq:step4_stima03}
\iint_{Q_{R}(x_{0},t_{0})}(1+\mathcal{U}^{p+2M})\,dx\,dt\,&\leq\,2^{n+1}+2^{p}\iint_{Q_{R}(x_{0},t_{0})\cap\{\mathcal{U}\,\geq\,1\}}\left(\mathcal{U}-\frac{1}{2}\right)_{+}^{p}\mathcal{U}^{2M}\,dx\,dt\nonumber\\
&\leq\,2^{n+1}+2^{p}\iint_{Q_{R}(x_{0},t_{0})}\left(\mathcal{U}-\frac{1}{2}\right)_{+}^{p}\mathcal{U}^{2M}\,dx\,dt\,.
\end{align}Joining (\ref{eq:st4_stima2}) and \eqref{eq:step4_stima03}, we then
find\begin{align}\label{eq:step4_final}
&\iint_{Q_{r}(x_{0},t_{0})}\mathcal{U}^{2M\,+\,\frac{4(M+1)}{n}}\left(\mathcal{U}-\frac{1}{2}\right)_{+}^{p}dx\,dt\nonumber\\
&\,\,\,\,\,\,\,\leq\left[C\,\frac{M^{6}}{(R-r)^{p}}\left(1+\iint_{Q_{R}(x_{0},t_{0})}\mathcal{U}^{2M}\left(\mathcal{U}-\frac{1}{2}\right)_{+}^{p}dx\,dt\right)\right]^{\frac{2}{n}+1}.
\end{align}\\
\textbf{Step 5: the local $L^{\infty}$ estimate on $Du$ in the case
$p>2$.} We now take $M=M_{j}=2^{j+1}-1$ with $j\in\mathbb{N}$.
Then, we set 
\[
\gamma_{j}:=\,2\,M_{j}=\,2^{j+2}-2,\,\,\,\,\,\,\,\,\widehat{\gamma}_{j}:=\,2\,M_{j}+\,\frac{4(M_{j}+1)}{n}\,=\,2^{j+2}-2+\,\frac{4}{n}\,2^{j+1},\,\,\,\,\,\,\,\,\,\,\,\,\mathrm{for}\,\,j\in\mathbb{N},
\]
and 
\[
\tau_{j}:=\,\frac{\widehat{\gamma}_{j}-\gamma_{j}}{\widehat{\gamma}_{j}-\gamma_{j-1}}\,\,\frac{\gamma_{j-1}}{\gamma_{j}}\,,\,\,\,\,\,\,\,\,\,\,\,\,\mathrm{for}\,\,j\in\mathbb{N}\,\backslash\,\{0\}.
\]
We note that $\gamma_{j-1}<\gamma_{j}<\widehat{\gamma}_{j}$ and $\tau_{j}\in(0,1)$
is defined in such a way that 
\[
\frac{1}{\gamma_{j}}\,=\,\frac{\tau_{j}}{\gamma_{j-1}}\,+\,\frac{1-\tau_{j}}{\widehat{\gamma}_{j}}\,.
\]
In order to simplify our notation, we also introduce the absolutely
continuous measure
\begin{equation}
d\mu:=\left(\mathcal{U}-\frac{1}{2}\right)_{+}^{p}d\mathcal{L}^{n+1}\,,\label{eq:mu}
\end{equation}
where $\mathcal{L}^{n+1}$ denotes the $(n+1)$-dimensional Lebesgue
measure. Thus, estimate \eqref{eq:step4_final} can be rewritten as
follows: 
\[
\int_{Q_{r}(x_{0},t_{0})}\mathcal{U}^{\widehat{\gamma}_{j}}\,d\mu\,\leq\left[C\,\frac{M_{j}^{6}}{(R-r)^{p}}\left(1+\int_{Q_{R}(x_{0},t_{0})}\mathcal{U}^{\gamma_{j}}\,d\mu\right)\right]^{\frac{2}{n}+1}.
\]
By interpolation in Lebesgue spaces, we obtain 
\[
\int_{Q_{r}(x_{0},t_{0})}\mathcal{U}^{\gamma_{j}}\,d\mu\,\leq\left(\int_{Q_{r}(x_{0},t_{0})}\mathcal{U}^{\gamma_{j-1}}\,d\mu\right)^{\tau_{j}\,\frac{\gamma_{j}}{\gamma_{j-1}}}\left(\int_{Q_{r}(x_{0},t_{0})}\mathcal{U}^{\widehat{\gamma}_{j}}\,d\mu\right)^{(1-\tau_{j})\,\frac{\gamma_{j}}{\widehat{\gamma}_{j}}}.
\]
Now, a few elementary computations reveal that 
\[
\tau_{j}\,\frac{\gamma_{j}}{\gamma_{j-1}}\,=\,\frac{4}{n+4}\,\,\,\,\,\,\,\,\,\,\,\,\,\,\,\,\mathrm{and}\,\,\,\,\,\,\,\,\,\,\,\,\,\,\,\,(1-\tau_{j})\,\frac{\gamma_{j}}{\widehat{\gamma}_{j}}\,=\,\frac{n}{n+4}\,.
\]
Thus, the combination of the two previous inequalities leads to 
\[
\int_{Q_{r}(x_{0},t_{0})}\mathcal{U}^{\gamma_{j}}\,d\mu\,\leq\left(\int_{Q_{R}(x_{0},t_{0})}\mathcal{U}^{\gamma_{j-1}}\,d\mu\right)^{\frac{4}{n+4}}\left[C\,\frac{M_{j}^{6}}{(R-r)^{p}}\left(1+\int_{Q_{R}(x_{0},t_{0})}\mathcal{U}^{\gamma_{j}}\,d\mu\right)\right]^{\frac{n+2}{n+4}}.
\]
By Young's inequality, we get\begin{align}\label{eq:step5_initial}
\int_{Q_{r}(x_{0},t_{0})}\mathcal{U}^{\gamma_{j}}\,d\mu\,&\leq\,\frac{n+2}{n+4}\int_{Q_{R}(x_{0},t_{0})}\mathcal{U}^{\gamma_{j}}\,d\mu\nonumber\\
&\,\,\,\,\,\,\,+\,\frac{2}{n+4}\left(\frac{CM_{j}^{6}}{(R-r)^{p}}\right)^{\frac{n+2}{2}}\left(\int_{Q_{R}(x_{0},t_{0})}\mathcal{U}^{\gamma_{j-1}}\,d\mu\right)^{2}+\,\frac{n+2}{n+4}\,.
\end{align} We can now invoke Lemma \ref{lem:giusti} to absorb the term on the
right-hand side of \eqref{eq:step5_initial} involving $\mathcal{U}^{\gamma_{j}}$,
in a standard way. By using the definition of $M_{j}$ and the fact
that $R\leq1$, we get\begin{align*}
\int_{Q_{r}(x_{0},t_{0})}\mathcal{U}^{\gamma_{j}}\,d\mu\,&\leq\,C\,\frac{2^{3(n+2)j}}{(R-r)^{p\,\frac{n+2}{2}}}\left(\int_{Q_{R}(x_{0},t_{0})}\mathcal{U}^{\gamma_{j-1}}\,d\mu\right)^{2}+C\nonumber\\
&\leq\,C\,\frac{2^{3(n+2)j}}{(R-r)^{p\,\frac{n+2}{2}}}\left[1+\left(\int_{Q_{R}(x_{0},t_{0})}\mathcal{U}^{\gamma_{j-1}}\,d\mu\right)^{2}\right]\nonumber\\
&\leq\,C\,\frac{2^{3(n+2)j}}{(R-r)^{p\,\frac{n+2}{2}}}\left(1+\int_{Q_{R}(x_{0},t_{0})}\mathcal{U}^{\gamma_{j-1}}\,d\mu\right)^{2},
\end{align*}for some constant $C=C(n,p,\delta)>1$. By summing $1$ on both sides
of the previous estimate, and exploiting that 
\[
C\,\frac{2^{3(n+2)j}}{(R-r)^{p\,\frac{n+2}{2}}}\,>\,1\,\,\,\,\,\,\,\,\,\,\,\,\mathrm{for\,\,every\,\,}j\in\mathbb{N}\,\backslash\,\{0\},
\]
we obtain 
\begin{equation}
1+\int_{Q_{r}(x_{0},t_{0})}\mathcal{U}^{\gamma_{j}}\,d\mu\,\leq\,C\,\frac{2^{3(n+2)j}}{(R-r)^{p\,\frac{n+2}{2}}}\left[1+\left(1+\int_{Q_{R}(x_{0},t_{0})}\mathcal{U}^{\gamma_{j-1}}\,d\mu\right)^{2}\right].\label{eq:iteranda}
\end{equation}
Now we want to iterate the above estimate on a sequence of shrinking
parabolic cylinders. To this end, we consider the decreasing sequence
\[
R_{j}:=\,r+\,\frac{R-r}{2^{j-1}}\,,\,\,\,\,\,\,\,\,\,\,\,\,j\in\mathbb{N}\,\backslash\,\{0\},
\]
and apply (\ref{eq:iteranda}) with $R_{j+1}<R_{j}$ in place of $r<R$.
To simplify our notation, we define
\begin{equation}
Y_{j}:=\,1+\int_{Q_{R_{j}}(x_{0},t_{0})}\mathcal{U}^{\gamma_{j-1}}\,d\mu\,,\,\,\,\,\,\,\,\,\,\,\,\,j\in\mathbb{N}\,\backslash\,\{0\},\label{eq:Y_j}
\end{equation}
and 
\begin{equation}
\vartheta:=\,\frac{n+2}{2}\,.\label{eq:TETA}
\end{equation}
Using the fact that $Y_{j}\geq1$ for every $j\in\mathbb{N}\,\backslash\,\{0\}$,
and up to redefining the constant $C>1$, from (\ref{eq:iteranda})
we get 
\[
Y_{j+1}\,\leq\,C\,2^{3\,p\,(n+2)j}\,(R-r)^{-\,\vartheta p}\,Y_{j}^{2}\,,
\]
for any $j\in\mathbb{N}\,\backslash\,\{0\}$. By iterating the previous
inequality starting from $j=1$, we obtain for every $k\in\mathbb{N}\,\backslash\,\{0\}$\begin{align*}
Y_{k+1}\,&\leq\,Y_{1}^{2^{k}}\prod_{j=1}^{k}\left[C\,8^{p\,(n+2)j}\,(R-r)^{-\,\vartheta p}\right]^{2^{k\,-\,j}}\\
&=\,Y_{1}^{2^{k}}\left[C\,(R-r)^{-\,\vartheta p}\right]^{\sum_{j=1}^{k}2^{k\,-\,j}}8^{p\,(n+2)\,\sum_{j=1}^{k}j\,2^{k\,-\,j}}.
\end{align*}Now we observe that 
\[
\sum_{j=1}^{k}2^{k-j}\,=\,2^{k}\,\sum_{j=1}^{k}\left(\frac{1}{2}\right)^{j}<\,2^{k}
\]
and 
\[
\sum_{j=1}^{k}j\,2^{k-j}\,=\,2^{k}\,\sum_{j=1}^{k}\,j\left(\frac{1}{2}\right)^{j}\leq\,2^{k}\,\frac{\left(\frac{1}{2}\right)}{\left(1-\frac{1}{2}\right)^{2}}\,=\,2^{k+1}.
\]
Therefore, for every $k\in\mathbb{N}\,\backslash\,\{0\}$, we have
\[
Y_{k+1}\,\leq\left[C\,64^{p\,(n+2)}\,(R-r)^{-\,\vartheta p}\,Y_{1}\right]^{2^{k}}.
\]
Thus, by redefining the constant $C>1$ again and recalling the definition
of $Y_{j}$ in (\ref{eq:Y_j}), we obtain
\[
\int_{Q_{R_{k+1}}(x_{0},t_{0})}\mathcal{U}^{\gamma_{k}}\,d\mu\,<\,Y_{k+1}\,\leq\left[\frac{C}{(R-r)^{\vartheta p}}\left(1+\int_{Q_{R}(x_{0},t_{0})}\mathcal{U}^{2}\,d\mu\right)\right]^{2^{k}},
\]
for any $k\in\mathbb{N}\,\backslash\,\{0\}$. Taking both sides of
the previous inequality to the power $\gamma_{k}^{-1}$ and then letting
$k\to\infty$, we find 
\[
\Vert\mathcal{U}\Vert_{L^{\infty}(Q_{r}(x_{0},t_{0}),\,d\mu)}\,\leq\,\frac{C}{(R-r)^{p\,\frac{n+2}{8}}}\left(1+\int_{Q_{R}(x_{0},t_{0})}\mathcal{U}^{2}\,d\mu\right)^{\frac{1}{4}}.
\]
Here, we have also used the definition of $\vartheta$ in (\ref{eq:TETA})
and the fact that $\gamma_{k}\sim2^{k+2}$ as $k$ tends to $\infty$.
Exploiting once again the condition $R\leq1$ and recalling the definition
of $d\mu$ in (\ref{eq:mu}), we deduce from the previous estimate
that\begin{align*}
\Vert\mathcal{U}\Vert_{L^{\infty}(Q_{r}(x_{0},t_{0}),\,d\mu)}\,&\leq\,\frac{C}{(R-r)^{p\,\frac{n+2}{4}}}\left(1+\iint_{Q_{R}(x_{0},t_{0})}\mathcal{U}^{2}\left(\mathcal{U}-\frac{1}{2}\right)_{+}^{p}dx\,dt\right)^{\frac{1}{4}}\\
&\leq\,\frac{C}{(R-r)^{p\,\frac{n+2}{4}}}\left(1+\iint_{Q_{R}(x_{0},t_{0})}\mathcal{U}^{p+2}\,dx\,dt\right)^{\frac{1}{4}}\\
&\leq\,\frac{C}{(R-r)^{p\,\frac{n+2}{4}}}\left(1+\iint_{Q_{R}(x_{0},t_{0})}\vert Du\vert^{p+2}\,dx\,dt\right)^{\frac{1}{4}},
\end{align*} where, in the last line, we have applied the inequalities $\mathcal{U}\leq\frac{1}{2\delta}\vert Du\vert\leq\vert Du\vert$.
Recalling the definition of $d\mu$ again, using the above estimate
and taking into account that $R\leq1<C$, we also get \begin{align*}
\Vert\mathcal{U}\Vert_{L^{\infty}(Q_{r}(x_{0},t_{0}))}\,&\leq\,\max\left\{ \Vert\mathcal{U}\Vert_{L^{\infty}(Q_{r}(x_{0},t_{0}),\,d\mu)}\,,\underset{Q_{r}(x_{0},t_{0})\cap\left\{ \mathcal{U}\,\le\,\frac{1}{2}\right\} }{\sup}\mathcal{U}\right\}\\
&\leq\,\max\left\{ \Vert\mathcal{U}\Vert_{L^{\infty}(Q_{r}(x_{0},t_{0}),\,d\mu)}\,,\frac{1}{2}\right\}\\
&\leq\,\frac{1}{2}\,+\,\Vert\mathcal{U}\Vert_{L^{\infty}(Q_{r}(x_{0},t_{0}),\,d\mu)}\\
&\leq\,\frac{C}{(R-r)^{p\,\frac{n+2}{4}}}\left(1+\iint_{Q_{R}(x_{0},t_{0})}\vert Du\vert^{p+2}\,dx\,dt\right)^{\frac{1}{4}}.
\end{align*}We now apply the second inequality in (\ref{eq:controlDu}) together
with the previous estimate. This yields \begin{align}\label{eq:step5_key}
\Vert Du\Vert_{L^{\infty}(Q_{r}(x_{0},t_{0}))}\,&\leq\,\sqrt{n}\,2\delta\,\Vert\mathcal{U}\Vert_{L^{\infty}(Q_{r}(x_{0},t_{0}))}\nonumber\\
&\leq\,\frac{C}{(R-r)^{p\,\frac{n+2}{4}}}\left[1+\left(\iint_{Q_{R}(x_{0},t_{0})}\vert Du\vert^{p+2}\,dx\,dt\right)^{\frac{1}{4}}\right].
\end{align}Finally, in order to remove the dependence on the $L^{p+2}$ norm
of the gradient, we use a standard interpolation argument. We write
\[
\left(\iint_{Q_{R}(x_{0},t_{0})}\vert Du\vert^{p+2}\,dx\,dt\right)^{\frac{1}{4}}\leq\,\Vert Du\Vert_{L^{\infty}(Q_{R}(x_{0},t_{0}))}^{\frac{1}{2}}\left(\iint_{Q_{R}(x_{0},t_{0})}\vert Du\vert^{p}\,dx\,dt\right)^{\frac{1}{4}}.
\]
Inserting this estimate into \eqref{eq:step5_key} and applying Young's
inequality, we obtain\begin{align*}
\Vert Du\Vert_{L^{\infty}(Q_{r}(x_{0},t_{0}))}\,&\leq\,\frac{C}{(R-r)^{p\,\frac{n+2}{4}}}\,+\,\frac{C}{(R-r)^{p\,\frac{n+2}{4}}}\,\Vert Du\Vert_{L^{\infty}(Q_{R}(x_{0},t_{0}))}^{\frac{1}{2}}\left(\iint_{Q_{R}(x_{0},t_{0})}\vert Du\vert^{p}\,dx\,dt\right)^{\frac{1}{4}}\\
&\leq\,\frac{1}{2}\,\Vert Du\Vert_{L^{\infty}(Q_{R}(x_{0},t_{0}))}\,+\,\frac{C}{(R-r)^{p\,\frac{n+2}{2}}}\left[1+\left(\iint_{Q_{R}(x_{0},t_{0})}\vert Du\vert^{p}\,dx\,dt\right)^{\frac{1}{2}}\right],
\end{align*} where, in the last line, we have also used that $R\leq1$. By Lemma
\ref{lem:giusti} again, we get
\[
\Vert Du\Vert_{L^{\infty}(Q_{r}(x_{0},t_{0}))}\,\leq\,\frac{C}{(R-r)^{p\,\frac{n+2}{2}}}\left[1+\left(\iint_{Q_{R}(x_{0},t_{0})}\vert Du\vert^{p}\,dx\,dt\right)^{\frac{1}{2}}\right].
\]
This concludes the proof for the case $p>2$.\\
\\
\textbf{Step 6: the local $L^{\infty}$ estimate on $Du$ in the case
$p=2$.} We now detail the modifications of the above proof to obtain
the local estimate (\ref{eq:uniform}) in the case $p=2$.

\noindent $\hspace*{1em}$Using the definition of $F$ in (\ref{eq:F_quad}),
the fact that $\widetilde{g}_{i,\varepsilon}''\geq\mathds{1}_{\{\vert\tau\vert\,\ge\,\delta_{i}\,+\,\varepsilon\}}$
and (\ref{eq:ell_F_quad}), we obtain that 
\[
\sum_{i\,\in\,E^{+}}\mathds{1}_{\{\vert\tau\vert\,\ge\,\delta_{i}\,+\,\varepsilon\}}(\xi_{i})\,\zeta_{i}^{2}\,+\sum_{i\,\in\,E^{-}}\zeta_{i}^{2}\,\leq\,\langle D^{2}F(\xi)\,\zeta,\zeta\rangle\,\leq\,2\,\vert\zeta\vert^{2}\,\,\,\,\,\,\,\,\,\,\,\,\mathrm{for\,\,every\,\,}\xi,\zeta\in\mathbb{R}^{n},
\]
where 
\begin{equation}
E^{+}:=\{i\in\{1,\ldots,n\}:\delta_{i}>0\}\,\,\,\,\,\,\,\,\,\,\,\,\mathrm{and}\,\,\,\,\,\,\,\,\,\,\,\,E^{-}:=\{i\in\{1,\ldots,n\}:\delta_{i}=0\}\,.\label{eq:EE}
\end{equation}
Inserting the previous estimates into \eqref{eq:step3_stima01}, one
gets\begin{align*}
&\chi(\tau)\,\sum_{\ell=0}^{\ell_{0}-1}\int_{\Omega\times\{\tau\}}\vert u_{x_{j}}\vert^{2s_{\ell}}\,\vert u_{x_{k}}\vert^{2m_{\ell}}\,\eta^{2}\,dx\\
&\,\,\,\,\,\,\,+\sum_{i\,\in\,E^{+}}\iint_{\Omega\times(T_{0},\tau)}\mathds{1}_{\{\vert u_{x_{i}}\vert\,\ge\,\delta_{i}\,+\,\varepsilon\}}\,u_{x_{i}x_{j}}^{2}\,\vert u_{x_{k}}\vert^{2M}\,\chi\,\eta^{2}\,dx\,dt\\
&\,\,\,\,\,\,\,+\sum_{i\,\in\,E^{-}}\iint_{\Omega\times(T_{0},\tau)}u_{x_{i}x_{j}}^{2}\,\vert u_{x_{k}}\vert^{2M}\,\chi\,\eta^{2}\,dx\,dt\\
&\,\,\,\,\,\,\,\,\,\,\,\,\,\,\leq\,522\,M^{4}\iint_{\Omega\times(T_{0},\tau)}[(\partial_{t}\chi)\,\eta^{2}+\chi\,\vert D\eta\vert^{2}]\,[\vert u_{x_{j}}\vert^{2(M+1)}+\vert u_{x_{k}}\vert^{2(M+1)}]\,dx\,dt\,.
\end{align*}Now we consider the last two terms on the left-hand side. By keeping
in the sums only the term with $i=k$ and dropping the others, we
obtain\begin{align*}
&\chi(\tau)\,\sum_{\ell=0}^{\ell_{0}-1}\int_{\Omega\times\{\tau\}}\vert u_{x_{j}}\vert^{2s_{\ell}}\,\vert u_{x_{k}}\vert^{2m_{\ell}}\,\eta^{2}\,dx\,+\iint_{\Omega\times(T_{0},\tau)}\mathds{1}_{\{\vert u_{x_{k}}\vert\,\ge\,\delta_{k}\,+\,\varepsilon\}}\,u_{x_{k}x_{j}}^{2}\,\vert u_{x_{k}}\vert^{2M}\,\chi\,\eta^{2}\,dx\,dt\\
&\,\,\,\,\,\,\,\leq\,522\,M^{4}\iint_{\Omega\times(T_{0},\tau)}[(\partial_{t}\chi)\,\eta^{2}+\chi\,\vert D\eta\vert^{2}]\,[\vert u_{x_{j}}\vert^{2(M+1)}+\vert u_{x_{k}}\vert^{2(M+1)}]\,dx\,dt\,.
\end{align*}Note that, by Young's inequality, one has\begin{align*}
&\iint_{\Omega\times(T_{0},\tau)}\left|[(\vert u_{x_{k}}\vert-\delta_{k}-\varepsilon)_{+}\,\vert u_{x_{k}}\vert^{M}]_{x_{j}}\right|^{2}\chi\,\eta^{2}\,dx\,dt\\
&\,\,\,\,\,\,\,\leq\,2\iint_{\Omega\times(T_{0},\tau)}\left|[(\vert u_{x_{k}}\vert-\delta_{k}-\varepsilon)_{+}]_{x_{j}}\right|^{2}\vert u_{x_{k}}\vert^{2M}\,\chi\,\eta^{2}\,dx\,dt\\
&\,\,\,\,\,\,\,\,\,\,\,\,\,\,+\,2\iint_{\Omega\times(T_{0},\tau)}(\vert u_{x_{k}}\vert-\delta_{k}-\varepsilon)_{+}^{2}\left|[\vert u_{x_{k}}\vert^{M}]_{x_{j}}\right|^{2}\chi\,\eta^{2}\,dx\,dt\\
&\,\,\,\,\,\,\,\leq\,4M^{2}\iint_{\Omega\times(T_{0},\tau)}\mathds{1}_{\{\vert u_{x_{k}}\vert\,\ge\,\delta_{k}\,+\,\varepsilon\}}\,u_{x_{k}x_{j}}^{2}\vert u_{x_{k}}\vert^{2M}\,\chi\,\eta^{2}\,dx\,dt\,,
\end{align*}where, in the last line, we have used that $M\geq1$. Then, combining
the two previous estimates and summing over $j\in\{1,\ldots,n\}$
the resulting inequality, we get\begin{align*}
&\chi(\tau)\,\sum_{\ell=0}^{\ell_{0}-1}\int_{\Omega\times\{\tau\}}\sum_{j=1}^{n}\,\vert u_{x_{j}}\vert^{2s_{\ell}}\,\vert u_{x_{k}}\vert^{2m_{\ell}}\,\eta^{2}\,dx\,+\,\frac{1}{4M^{2}}\iint_{\Omega\times(T_{0},\tau)}\left|D[(\vert u_{x_{k}}\vert-\delta_{k}-\varepsilon)_{+}\,\vert u_{x_{k}}\vert^{M}]\right|^{2}\chi\,\eta^{2}\,dx\,dt\\
&\,\,\,\,\,\,\,\leq\,522\,M^{4}\iint_{\Omega\times(T_{0},\tau)}[(\partial_{t}\chi)\,\eta^{2}+\chi\,\vert D\eta\vert^{2}]\left[\sum_{j=1}^{n}\,\vert u_{x_{j}}\vert^{2M+2}+n\,\vert u_{x_{k}}\vert^{2M+2}\right]dx\,dt\,.
\end{align*}We now use that $4M^{2}\geq4$ and add the term 
\[
\iint_{\Omega\times(T_{0},\tau)}(\vert u_{x_{k}}\vert-\delta_{k}-\varepsilon)_{+}^{2}\,\vert u_{x_{k}}\vert^{2M}\,\chi\,\vert D\eta\vert^{2}\,dx\,dt
\]
to both sides of the preceding inequality. With some algebraic manipulations,
this gives\begin{align*}
&\chi(\tau)\,\sum_{\ell=0}^{\ell_{0}-1}\int_{\Omega\times\{\tau\}}\sum_{j=1}^{n}\,\vert u_{x_{j}}\vert^{2s_{\ell}}\,\vert u_{x_{k}}\vert^{2m_{\ell}}\,\eta^{2}\,dx\,+\iint_{\Omega\times(T_{0},\tau)}\left|D[(\vert u_{x_{k}}\vert-\delta_{k}-\varepsilon)_{+}\,\vert u_{x_{k}}\vert^{M}\,\eta]\right|^{2}\chi\,dx\,dt\\
&\,\,\,\,\,\,\,\leq\,c_{0}\,M^{6}\iint_{\Omega\times(T_{0},\tau)}[(\partial_{t}\chi)\,\eta^{2}+\chi\,\vert D\eta\vert^{2}]\left[\sum_{j=1}^{n}\,\vert u_{x_{j}}\vert^{2M+2}\,+\,\vert u_{x_{k}}\vert^{2M+2}\right]dx\,dt\\
&\,\,\,\,\,\,\,\,\,\,\,\,\,\,+\,2\iint_{\Omega\times(T_{0},\tau)}(\vert u_{x_{k}}\vert-\delta_{k}-\varepsilon)_{+}^{2}\,\vert u_{x_{k}}\vert^{2M}\,\chi\,\vert D\eta\vert^{2}\,dx\,dt\\
&\,\,\,\,\,\,\,\leq\,c\,M^{6}\iint_{\Omega\times(T_{0},\tau)}[(\partial_{t}\chi)\,\eta^{2}+\chi\,\vert D\eta\vert^{2}]\left[\sum_{j=1}^{n}\,\vert u_{x_{j}}\vert^{2M+2}\,+\,\vert u_{x_{k}}\vert^{2M+2}\right]dx\,dt\,,
\end{align*}where $c_{0}$ and $c$ are positive constants depending only on $n$.
By using the Sobolev inequality in the spatial variable for the second
term on the left-hand side, we obtain\begin{align*}
&\chi(\tau)\,\sum_{\ell=0}^{\ell_{0}-1}\int_{\Omega\times\{\tau\}}\sum_{j=1}^{n}\,\vert u_{x_{j}}\vert^{2s_{\ell}}\,\vert u_{x_{k}}\vert^{2m_{\ell}}\,\eta^{2}\,dx\,+\int_{T_{0}}^{\tau}\chi\left(\int_{\Omega}(\vert u_{x_{k}}\vert-\delta_{k}-\varepsilon)_{+}^{2^{*}}\,\vert u_{x_{k}}\vert^{2^{*}M}\,\eta^{2^{*}}dx\right)^{\frac{2}{2^{*}}}dt\\
&\,\,\,\,\,\,\,\leq\,c\,M^{6}\iint_{\Omega\times(T_{0},\tau)}[(\partial_{t}\chi)\,\eta^{2}+\chi\,\vert D\eta\vert^{2}]\left[\sum_{j=1}^{n}\,\vert u_{x_{j}}\vert^{2M+2}\,+\,\vert u_{x_{k}}\vert^{2M+2}\right]dx\,dt\,.
\end{align*}At this stage, we sum over $k\in\{1,\ldots,n\}$ and apply Minkowski's
inequality to the second term on the left-hand side. This yields\begin{align}\label{eq:step6_stima02}
&\chi(\tau)\,\sum_{\ell=0}^{\ell_{0}-1}\int_{\Omega\times\{\tau\}}\sum_{j=1}^{n}\,\vert u_{x_{j}}\vert^{2s_{\ell}}\,\sum_{k=1}^{n}\,\vert u_{x_{k}}\vert^{2m_{\ell}}\,\eta^{2}\,dx\nonumber\\
&\,\,\,\,+\int_{T_{0}}^{\tau}\chi\left(\int_{\Omega}\left|\sum_{k=1}^{n}\,(\vert u_{x_{k}}\vert-\delta_{k}-\varepsilon)_{+}^{2}\,\vert u_{x_{k}}\vert^{2M}\right|^{\frac{2^{*}}{2}}\eta^{2^{*}}dx\right)^{\frac{2}{2^{*}}}dt\nonumber\\
&\,\,\,\,\,\,\,\,\,\,\,\leq\,c\,M^{6}\iint_{\Omega\times(T_{0},\tau)}[(\partial_{t}\chi)\,\eta^{2}+\chi\,\vert D\eta\vert^{2}]\left[\sum_{j=1}^{n}\,\vert u_{x_{j}}\vert^{2M+2}+\sum_{k=1}^{n}\,\vert u_{x_{k}}\vert^{2M+2}\right]dx\,dt.
\end{align}Similarly to what has been done in steps 3 and 5, we now introduce
the auxiliary function 
\[
\mathcal{V}(x,t):=\,\frac{1}{2(\delta+1)}\,\,\underset{1\,\le\,i\,\leq\,n}{\max}\vert u_{x_{i}}(x,t)\vert
\]
and the absolutely continuous measure
\[
d\sigma:=\left(\mathcal{V}-\frac{1}{2}\right)_{+}^{2}d\mathcal{L}^{n+1}\,.
\]
A few elementary computations reveal that 
\[
\sum_{k=1}^{n}(\vert u_{x_{k}}\vert-\delta_{k}-\varepsilon)_{+}^{2}\,\vert u_{x_{k}}\vert^{2M}\,\ge\,(2\delta+2)^{2M+2}\left(\mathcal{V}-\frac{1}{2}\right)_{+}^{2}\,\mathcal{V}^{2M}
\]
and
\[
[2(\delta+1)\,\mathcal{V}]^{q}\,\le\,\sum_{k=1}^{n}\,\vert u_{x_{k}}\vert^{q}\,\leq\,n\,[2(\delta+1)\,\mathcal{V}]^{q}\,\,\,\,\,\,\,\,\,\,\,\,\mathrm{for\,\,every}\,\,q\geq0\,.
\]
In particular, for $q=2$ we obtain 
\[
2(\delta+1)\,\mathcal{V}\,\leq\,\vert Du\vert\,\leq\,\sqrt{n}\,2(\delta+1)\,\mathcal{V}\,.
\]
Inserting these estimates into \eqref{eq:step6_stima02}, using that
$\ell_{0}\geq1$ and recalling that $s_{\ell}+m_{\ell}=M+1$ for every
$\ell\in\{0,\ldots,\ell_{0}-1\}$, we get\begin{align*}
&\chi(\tau)\int_{\Omega\times\{\tau\}}\mathcal{V}^{2M+2}\,\eta^{2}\,dx\,+\int_{T_{0}}^{\tau}\chi\left(\int_{\Omega}\left(\mathcal{V}-\frac{1}{2}\right)_{+}^{2^{*}}\,\mathcal{V}^{2^{*}M}\,\eta^{2^{*}}dx\right)^{\frac{2}{2^{*}}}dt\nonumber\\
&\,\,\,\,\,\,\,\leq\,c\,M^{6}\iint_{\Omega\times(T_{0},\tau)}[(\partial_{t}\chi)\,\eta^{2}+\chi\,\vert D\eta\vert^{2}]\,\mathcal{V}^{2M+2}\,dx\,dt\,.
\end{align*} Starting from this estimate and proceeding exactly as in steps 4
and 5, but using $\mathcal{V}$, $\delta+1$ and $d\sigma$ in place
of $\mathcal{U}$, $\delta$ and $d\mu$, respectively, we reach the
desired conclusion for $p=2$.\end{proof}

\noindent \begin{brem}A careful inspection of the previous proof
reveals that the exponent $\vartheta$ in (\ref{eq:uniform}) can
be taken to be
\begin{equation}
\vartheta=\begin{cases}
\begin{array}{cc}
{\displaystyle \frac{n+2}{2}} & \,\,\mathrm{if}\,\,n\geq3,\vspace{1mm}\\
\mathrm{any\,\,number}>2 & \,\,\mathrm{if}\,\,n=2.
\end{array}\end{cases}\label{eq:theta}
\end{equation}
In the case $n=2$, the constant $C$ in (\ref{eq:uniform}) blows
up as $\vartheta\to2$.\end{brem}

\section{Uniform energy estimates for a regularized problem\label{sec:energy} }

$\hspace*{1em}$Let us fix an open set $\Omega'\Subset\Omega$ and
a subinterval $J:=(\tau_{0},\tau_{1})\Subset I$. Assume that $u\in L_{loc}^{p}(I;W_{loc}^{1,p}(\Omega))$
is a local weak solution of equation (\ref{eq:equation}). In light
of (\ref{eq:additional_properties-1}), we have 
\[
\partial_{t}u\in L^{p'}(J;W^{-1,p'}(\Omega'))\,\,\,\,\,\,\,\,\,\,\,\,\,\,\,\mathrm{and}\,\,\,\,\,\,\,\,\,\,\,\,\,\,\,u\in C^{0}(\overline{J};L^{2}(\Omega'))\,.
\]
Now, for any fixed $\varepsilon\in(0,\min\,\{1,\inf\Delta^{+}\})$,
we consider the approximating Cauchy-Dirichlet problem parametrized
by $\varepsilon$
\begin{equation}
\begin{cases}
\begin{array}{cc}
\partial_{t}v\,=\,\mathrm{div}\,[D_{\xi}F_{\varepsilon}(Dv)] & \mathrm{in\,\,\,}\Omega'\times J,\\
v=u & \,\,\,\,\mathrm{on\,\,\,}\partial\Omega'\times J,\\
v(\cdot,\tau_{0})=u(\cdot,\tau_{0}) & \mathrm{in\,\,\,}\Omega'\,.\,\,\,\,\,\,\,\,
\end{array}\end{cases}\label{eq:regularized}
\end{equation}
By \cite[Proposition 4.1, Chapter III]{Show}, this problem admits
a unique weak solution $u_{\varepsilon}\in L^{p}(J;W^{1,p}(\Omega'))$
such that 
\[
\partial_{t}u_{\varepsilon}\in L^{p'}(J;W^{-1,p'}(\Omega'))\,,\,\,\,\mathrm{and\,\,thus}\,\,\,u_{\varepsilon}\in C^{0}(\overline{J};L^{2}(\Omega'))\,.
\]
The condition $u_{\varepsilon}=u$ on the lateral boundary $\partial\Omega'\times J$
is understood in the sense that 
\[
u_{\varepsilon}-u\,\in\,L^{p}(J;W_{0}^{1,p}(\Omega'))\,,
\]
while the initial condition $u_{\varepsilon}(\cdot,\tau_{0})=u(\cdot,\tau_{0})$
in $\Omega'$ is taken in the usual $L^{2}$-sense, which is feasible
due to the continuity properties of both $u_{\varepsilon}$ and $u$.\\
$\hspace*{1em}$The penultimate step in the proof of Theorem \ref{thm:main}
consists in establishing the uniform energy estimates, as well as
the strong convergence results, stated in Propositions \ref{prop:unif_en}
and \ref{prop:unif_en-1} below. The need to distinguish the cases
$p>2$ and $p=2$ arises from the fact that, for $p=2$, the regularizing
function $F_{\varepsilon}$ is defined by (\ref{eq:F_quad}) rather
than (\ref{eq:F_eps}).
\begin{prop}[\textbf{Uniform energy estimate for $p>2$}]
\label{prop:unif_en}Let $n\geq2$, $p>2$ and $\varepsilon\in(0,\min\,\{1,\inf\Delta^{+}\})$.
Moreover, let $u\in L_{loc}^{p}(I;W_{loc}^{1,p}(\Omega))$ be a local
weak solution of $(\ref{eq:equation})$ and assume that $u_{\varepsilon}\in L^{p}(J;W^{1,p}(\Omega'))$
is the unique weak solution of problem $(\ref{eq:regularized})$.
Then, the estimate\begin{align}\label{eq:EnEst}
&\int_{\Omega '\times\{\tau_{1}\}}\vert u_{\varepsilon}-u\vert^{2}\,dx\,+\,\sum_{i=1}^{n}\iint_{\Omega '\times J}\left|H_{\delta_{i}}(\partial_{x_{i}}u_{\varepsilon})-H_{\delta_{i}}(\partial_{x_{i}}u)\right|^{2}dx\,dt\,+\,\varepsilon\iint_{\Omega '\times J}\vert Du_{\varepsilon}\vert^{p}\,dx\,dt\nonumber\\
&\,\,\,\,\,\,\,\leq\,c\,\varepsilon\left(\vert\Omega '\times J\vert+\iint_{\Omega '\times J}\vert Du\vert^{p}\,dx\,dt\right)
\end{align}holds for some positive constant $c$ depending only on $p$. In particular,
this estimate implies that 
\begin{equation}
\iint_{\Omega'\times J}\vert Du_{\varepsilon}\vert^{p}\,dx\,dt\,\leq\,c\left(\vert\Omega'\times J\vert+\iint_{\Omega'\times J}\vert Du\vert^{p}\,dx\,dt\right)\label{eq:uni}
\end{equation}
and
\begin{equation}
H_{\delta_{j}}(\partial_{x_{j}}u_{\varepsilon})\rightarrow H_{\delta_{j}}(\partial_{x_{j}}u)\,\,\,\,\,\,\,\,\mathit{strongly\,\,\,in}\,\,\,L^{2}(\Omega'\times J)\,\,\,\,\,\mathit{as}\,\,\,\varepsilon\to0\,,\label{eq:convergence}
\end{equation}
for each $j\in\{1,\ldots,n\}$. 
\end{prop}

\noindent \begin{proof}[\bfseries{Proof}]The function $u_{\varepsilon}$
verifies 
\begin{equation}
\iint_{\Omega'\times J}u_{\varepsilon}\,\partial_{t}\varphi\,dx\,dt\,-\,\iint_{\Omega'\times J}\langle D_{\xi}F_{\varepsilon}(Du_{\varepsilon}),D\varphi\rangle\,dx\,dt\,=\,0\,,\label{eq:weak_appro}
\end{equation}
for every $\varphi\in C_{0}^{\infty}(\Omega'\times J).$ Integrating
by parts in (\ref{eq:weak_appro}) yields 
\[
\int_{J}(\partial_{t}u_{\varepsilon},\varphi)_{(W^{-1,p'},W_{0}^{1,p})}\,dt\,+\,\iint_{\Omega'\times J}\langle D_{\xi}F_{\varepsilon}(Du_{\varepsilon}),D\varphi\rangle\,dx\,dt\,=\,0\,.
\]
By density, the above identity also holds for every $\varphi\in L^{p}(J;W_{0}^{1,p}(\Omega'))$.
We then choose $\varphi=u_{\varepsilon}-u$, which gives
\[
\int_{J}(\partial_{t}u_{\varepsilon},u_{\varepsilon}-u)_{(W^{-1,p'},W_{0}^{1,p})}\,dt\,+\,\iint_{\Omega'\times J}\langle D_{\xi}F_{\varepsilon}(Du_{\varepsilon}),Du_{\varepsilon}-Du\rangle\,dx\,dt\,=\,0\,.
\]
By recalling the definition of $F_{\varepsilon}$ in (\ref{eq:F_eps}),
the previous integral identity can be rewritten as follows:\begin{align*}
&\int_{J}(\partial_{t}u_{\varepsilon},u_{\varepsilon}-u)_{(W^{-1,p'},W_{0}^{1,p})}\,dt\,+\iint_{\Omega '\times J}\langle D_{\xi}F_{0}(Du_{\varepsilon}),Du_{\varepsilon}-Du\rangle\,dx\,dt\\
&\,\,\,\,\,\,\,+\,\varepsilon\iint_{\Omega '\times J}\langle D_{\xi}G(Du_{\varepsilon}),Du_{\varepsilon}-Du\rangle\,dx\,dt\,=\,0\,.
\end{align*} Starting from (\ref{eq:loc_weak_sol}), we similarly have
\[
\int_{J}(\partial_{t}u,u_{\varepsilon}-u)_{(W^{-1,p'},W_{0}^{1,p})}\,dt\,+\iint_{\Omega'\times J}\langle D_{\xi}F_{0}(Du),Du_{\varepsilon}-Du\rangle\,dx\,dt\,=\,0\,.
\]
By subtracting the two identities above, we get\begin{align}\label{eq:subtraction}
&\int_{J}(\partial_{t}u_{\varepsilon}-\partial_{t}u,u_{\varepsilon}-u)_{(W^{-1,p'},W_{0}^{1,p})}\,dt\,+\iint_{\Omega '\times J}\langle D_{\xi}F_{0}(Du_{\varepsilon})-D_{\xi}F_{0}(Du),Du_{\varepsilon}-Du\rangle\,dx\,dt\nonumber\\
&\,\,\,\,\,\,\,+\,\varepsilon\iint_{\Omega '\times J}\langle D_{\xi}G(Du_{\varepsilon}),Du_{\varepsilon}-Du\rangle\,dx\,dt\,=\,0\,.
\end{align}The term involving the time derivatives can be rewritten as 
\[
\int_{J}(\partial_{t}u_{\varepsilon}-\partial_{t}u,u_{\varepsilon}-u)_{(W^{-1,p'},W_{0}^{1,p})}\,dt\,=\,\frac{1}{2}\int_{\Omega'\times\{\tau_{1}\}}\vert u_{\varepsilon}-u\vert^{2}\,dx\,.
\]
This follows from the fact that the map 
\[
t\,\mapsto\,\frac{1}{2}\int_{\Omega'}\vert u_{\varepsilon}(x,t)-u(x,t)\vert^{2}\,dx
\]
is absolutely continuous on $J$, with derivative given exactly by
\[
(\partial_{t}u_{\varepsilon}-\partial_{t}u,u_{\varepsilon}-u)_{(W^{-1,p'},W_{0}^{1,p})}\,\,\,\,\,\,\,\,\,\,\,\,\mathrm{for\,\,a.e.\,\,}t\in J,
\]
see \cite[Proposition 1.2, Chapter III]{Show}.\\
$\hspace*{1em}$For the second term in \eqref{eq:subtraction}, we
apply Lemma \ref{lem:Brasco}, which, for every $w,z\in\mathbb{R}^{n}$,
yields 
\[
\langle D_{\xi}F_{0}(w)-D_{\xi}F_{0}(z),w-z\rangle=\sum_{i=1}^{n}\left(J_{\delta_{i}}(w_{i})-J_{\delta_{i}}(z_{i})\right)(w_{i}-z_{i})\,\geq\,\frac{4}{p^{2}}\,\sum_{i=1}^{n}\left|H_{\delta_{i}}(w_{i})-H_{\delta_{i}}(z_{i})\right|^{2},
\]
where the functions $J_{\delta_{i}}$ and $H_{\delta_{i}}$ are defined
respectively by (\ref{eq:J}) and (\ref{eq:H}) with $\lambda=\delta_{i}$.
Using this pointwise estimate in \eqref{eq:subtraction} and recalling
the definition of $G$ in (\ref{eq:G}), we then get \begin{align}\label{eq:energy00}
&\frac{1}{2}\int_{\Omega '\times\{\tau_{1}\}}\vert u_{\varepsilon}-u\vert^{2}\,dx\,+\,\frac{4}{p^{2}}\,\sum_{i=1}^{n}\iint_{\Omega '\times J}\left|H_{\delta_{i}}(\partial_{x_{i}}u_{\varepsilon})-H_{\delta_{i}}(\partial_{x_{i}}u)\right|^{2}dx\,dt\nonumber\\
&\,\,\,\,\,\,\,+\,\varepsilon\iint_{\Omega '\times J}\langle(1+\vert Du_{\varepsilon}\vert^{2})^{\frac{p-2}{2}}Du_{\varepsilon},Du_{\varepsilon}-Du\rangle\,dx\,dt\,\leq\,0\,.
\end{align}By the Cauchy-Schwarz inequality and Young's inequality with $\beta>0$,
from \eqref{eq:energy00} we infer\begin{align}\label{eq:energy1}
&\int_{\Omega '\times\{\tau_{1}\}}\vert u_{\varepsilon}-u\vert^{2}\,dx\,+\,\frac{8}{p^{2}}\,\sum_{i=1}^{n}\iint_{\Omega '\times J}\left|H_{\delta_{i}}(\partial_{x_{i}}u_{\varepsilon})-H_{\delta_{i}}(\partial_{x_{i}}u)\right|^{2}dx\,dt\,+\,2\varepsilon\iint_{\Omega '\times J}\vert Du_{\varepsilon}\vert^{p}\,dx\,dt\nonumber\\
&\,\,\,\,\,\,\,\leq\,2\varepsilon\iint_{\Omega '\times J}(1+\vert Du_{\varepsilon}\vert^{2})^{\frac{p-1}{2}}\,\vert Du\vert\,dx\,dt\nonumber\\
&\,\,\,\,\,\,\,\leq\,\frac{2\,\varepsilon\beta^{p'}}{p'}\iint_{\Omega '\times J}(1+\vert Du_{\varepsilon}\vert^{2})^{\frac{p}{2}}\,dx\,dt\,+\,\frac{2\,\varepsilon}{p\beta^{p}}\iint_{\Omega '\times J}\vert Du\vert^{p}\,dx\,dt\nonumber\\
&\,\,\,\,\,\,\,\leq\,\frac{2^{\frac{p}{2}}\,\varepsilon\beta^{p'}}{p'}\iint_{\Omega '\times J}\vert Du_{\varepsilon}\vert^{p}\,dx\,dt\,+\,\frac{2^{\frac{p}{2}}\,\varepsilon\beta^{p'}}{p'}\,\vert\Omega '\times J\vert\,+\,\frac{2\,\varepsilon}{p\beta^{p}}\iint_{\Omega '\times J}\vert Du\vert^{p}\,dx\,dt\,.
\end{align}Upon choosing $\beta=\left(\frac{p'}{2^{p/2}}\right)^{\frac{1}{p'}}$
and absorbing the first integral on the right-hand side of \eqref{eq:energy1}
into the left-hand side, we arrive at estimate \eqref{eq:EnEst}.\\
$\hspace*{1em}$The uniform energy estimate (\ref{eq:uni}) follows
by discarding the first two terms on the left-hand side of \eqref{eq:EnEst}
and then dividing by $\varepsilon$. Similarly, by dropping the first
and third terms on the left-hand side of \eqref{eq:EnEst} and letting
$\varepsilon\rightarrow0$, we obtain the conclusion (\ref{eq:convergence}).\end{proof}

\noindent $\hspace*{1em}$In the case $p=2$, to obtain a result analogous
to Proposition \ref{prop:unif_en}, we need to introduce, for each
$i\in\{1,\ldots,n\}$, the auxiliary function 
\begin{equation}
K_{i,\varepsilon}(s):=\int_{0}^{s}\sqrt{\,\widetilde{g}_{i,\varepsilon}''(\tau)}\,d\tau\,,\,\,\,\,\,\,\,\,\,\,\,\,s\in\mathbb{R}\,,\label{eq:new_jolly}
\end{equation}
where $\widetilde{g}_{i,\varepsilon}$ denotes the convex $C^{2}$
map defined in (\ref{eq:g_eps}) and (\ref{eq:g_eps_tilde}). More
precisely, we have the following 
\begin{prop}[\textbf{Uniform energy estimate for $p=2$}]
\label{prop:unif_en-1}Let $n\geq2$, $p=2$ and $\varepsilon\in(0,\min\,\{1,\inf\Delta^{+}\})$.
Moreover, let $u\in L_{loc}^{2}(I;W_{loc}^{1,2}(\Omega))$ be a local
weak solution of $(\ref{eq:equation})$ and assume that $u_{\varepsilon}\in L^{2}(J;W^{1,2}(\Omega'))$
is the unique weak solution of problem $(\ref{eq:regularized})$,
where $F_{\varepsilon}$ is defined by $(\ref{eq:F_quad})$. Then,
the estimate\begin{align}\label{eq:EnEstN2}
&\int_{\Omega'\times\{\tau_{1}\}}\vert u_{\varepsilon}-u\vert^{2}\,dx\,+\sum_{i=1}^{n}\iint_{\Omega'\times J}\left|K_{i,\varepsilon}(\partial_{x_{i}}u_{\varepsilon})-K_{i,\varepsilon}(\partial_{x_{i}}u)\right|^{2}dx\,dt\,+\,\varepsilon\iint_{\Omega'\times J}\vert Du_{\varepsilon}-Du\vert^{2}\,dx\,dt\nonumber\\
&\,\,\,\,\,\,\,\leq\,c\,\varepsilon\left(\vert\Omega'\times J\vert\,+\iint_{\Omega'\times J}\vert Du\vert^{2}\,dx\,dt\right)
\end{align}holds for some positive constant $c$ depending only on $n$. In particular,
this estimate implies that 
\begin{equation}
\iint_{\Omega'\times J}\vert Du_{\varepsilon}\vert^{2}\,dx\,dt\,\le\,2(c+1)\left(\vert\Omega'\times J\vert\,+\iint_{\Omega'\times J}\vert Du\vert^{2}\,dx\,dt\right).\label{eq:uni_2}
\end{equation}
Furthermore, we have
\begin{equation}
K_{j,\varepsilon}(\partial_{x_{j}}u_{\varepsilon})\rightarrow H_{\delta_{j}}(\partial_{x_{j}}u)\,\,\,\,\,\,\,\,\mathit{strongly\,\,\,in}\,\,\,L^{2}(\Omega'\times J)\,\,\,\,\,\mathit{as}\,\,\,\varepsilon\to0\,,\label{eq:convergence-1}
\end{equation}
for each $j\in\{1,\ldots,n\}$. 
\end{prop}

\noindent \begin{proof}[\bfseries{Proof}]Arguing exactly as in the
first part of the preceding proof, we find that \begin{align*}
&\frac{1}{2}\int_{\Omega'\times\{\tau_{1}\}}\vert u_{\varepsilon}-u\vert^{2}\,dx\,+\iint_{\Omega'\times J}\langle D_{\xi}F_{\varepsilon}(Du_{\varepsilon})-D_{\xi}F_{\varepsilon}(Du),Du_{\varepsilon}-Du\rangle\,dx\,dt\\
&\,\,\,\,\,\,\,=\iint_{\Omega'\times J}\langle D_{\xi}F_{0}(Du)-D_{\xi}F_{\varepsilon}(Du),Du_{\varepsilon}-Du\rangle\,dx\,dt\,.
\end{align*} Recalling the definitions of $F_{\varepsilon}$, $\widetilde{g}_{i,\varepsilon}$
and $E^{+}$ in (\ref{eq:F_quad}), (\ref{eq:g_eps_tilde}) and (\ref{eq:EE}),
respectively, the above identity can be rewritten as follows:\begin{align*}
&\frac{1}{2}\int_{\Omega'\times\{\tau_{1}\}}\vert u_{\varepsilon}-u\vert^{2}\,dx\,+\sum_{i=1}^{n}\iint_{\Omega'\times J}\left[\widetilde{g}_{i,\varepsilon}'(\partial_{x_{i}}u_{\varepsilon})-\widetilde{g}_{i,\varepsilon}'(\partial_{x_{i}}u)\right]\left[\partial_{x_{i}}u_{\varepsilon}-\partial_{x_{i}}u\right]dx\,dt\\
&\,\,+\,\varepsilon\iint_{\Omega'\times J}\vert Du_{\varepsilon}-Du\vert^{2}\,dx\,dt\\
&\,\,\,\,\,\,\,\,\,=\sum_{i\,\in\,E^{+}}\iint_{\Omega'\times J}\left[g_{i}'(\partial_{x_{i}}u)-g_{i,\varepsilon}'(\partial_{x_{i}}u)\right]\left[\partial_{x_{i}}u_{\varepsilon}-\partial_{x_{i}}u\right]dx\,dt\,-\,\varepsilon\iint_{\Omega'\times J}\langle Du,Du_{\varepsilon}-Du\rangle\,dx\,dt\,.
\end{align*}We can estimate the second term on the left-hand side by applying
Lemma \ref{lem:banal} with $v=\widetilde{g}_{i,\varepsilon}$ and
recalling (\ref{eq:new_jolly}). Thus, we obtain\begin{align*}
&\sum_{i=1}^{n}\iint_{\Omega'\times J}\left[\widetilde{g}_{i,\varepsilon}'(\partial_{x_{i}}u_{\varepsilon})-\widetilde{g}_{i,\varepsilon}'(\partial_{x_{i}}u)\right]\left[\partial_{x_{i}}u_{\varepsilon}-\partial_{x_{i}}u\right]dx\,dt\\
&\,\,\,\,\,\,\,\geq\sum_{i=1}^{n}\iint_{\Omega'\times J}\left|K_{i,\varepsilon}(\partial_{x_{i}}u_{\varepsilon})-K_{i,\varepsilon}(\partial_{x_{i}}u)\right|^{2}dx\,dt\,.
\end{align*}Combining the two previous estimates and using the Cauchy-Schwarz
and Young's inequalities, we get\begin{align}\label{eq:energy2_00}
&\frac{1}{2}\int_{\Omega'\times\{\tau_{1}\}}\vert u_{\varepsilon}-u\vert^{2}\,dx\,+\sum_{i=1}^{n}\iint_{\Omega'\times J}\left|K_{i,\varepsilon}(\partial_{x_{i}}u_{\varepsilon})-K_{i,\varepsilon}(\partial_{x_{i}}u)\right|^{2}dx\,dt\,+\,\varepsilon\iint_{\Omega'\times J}\vert Du_{\varepsilon}-Du\vert^{2}\,dx\,dt\nonumber\\
&\,\,\,\,\,\,\,\leq\sum_{i\,\in\,E^{+}}\iint_{\Omega'\times J}\left|g_{i}'(\partial_{x_{i}}u)-g_{i,\varepsilon}'(\partial_{x_{i}}u)\right|\left|\partial_{x_{i}}u_{\varepsilon}-\partial_{x_{i}}u\right|dx\,dt\nonumber\\
&\,\,\,\,\,\,\,\,\,\,\,\,\,\,+\,\frac{\varepsilon}{2}\iint_{\Omega'\times J}\vert Du_{\varepsilon}-Du\vert^{2}\,dx\,dt\,+\,\frac{\varepsilon}{2}\iint_{\Omega'\times J}\vert Du\vert^{2}\,dx\,dt\,.
\end{align}To estimate the first term on the right-hand side, we apply the second
inequality of Lemma \ref{lem:lemmaA1} from the \hyperref[sec:appendice]{Appendix},
together with Young's inequality in the form
\[
ab\,\leq\,\frac{a^{2}}{\varepsilon}\,+\,\frac{\varepsilon\,b^{2}}{4}\,.
\]
Consequently, we obtain\begin{align}\label{eq:energy2_est1}
&\sum_{i\,\in\,E^{+}}\iint_{\Omega'\times J}\left|g_{i}'(\partial_{x_{i}}u)-g_{i,\varepsilon}'(\partial_{x_{i}}u)\right|\left|\partial_{x_{i}}u_{\varepsilon}-\partial_{x_{i}}u\right|dx\,dt\nonumber\\
&\,\,\,\,\,\,\,\leq\sum_{i\,\in\,E^{+}}\left[\frac{\varepsilon}{16}\,\vert\Omega'\times J\vert\,+\,\frac{\varepsilon}{4}\iint_{\Omega'\times J}\left|\partial_{x_{i}}u_{\varepsilon}-\partial_{x_{i}}u\right|^{2}dx\,dt\right]\nonumber\\
&\,\,\,\,\,\,\,\leq\,\frac{n\,\varepsilon}{16}\,\vert\Omega'\times J\vert\,+\,\frac{\varepsilon}{4}\iint_{\Omega'\times J}\vert Du_{\varepsilon}-Du\vert^{2}\,dx\,dt\,.
\end{align}Joining \eqref{eq:energy2_00} and \eqref{eq:energy2_est1}, we deduce
the estimate \eqref{eq:EnEstN2}. As immediate consequences of \eqref{eq:EnEstN2},
we get\begin{align*}
\iint_{\Omega'\times J}\vert Du_{\varepsilon}\vert^{2}\,dx\,dt\,&\leq\,2\iint_{\Omega'\times J}\vert Du_{\varepsilon}-Du\vert^{2}\,dx\,dt\,+\,2\iint_{\Omega'\times J}\vert Du\vert^{2}\,dx\,dt\\
&\le\,2(c+1)\left(\vert\Omega'\times J\vert\,+\iint_{\Omega'\times J}\vert Du\vert^{2}\,dx\,dt\right)\,,
\end{align*}and moreover 
\begin{equation}
{\displaystyle \big\Vert K_{j,\varepsilon}(\partial_{x_{j}}u_{\varepsilon})-K_{j,\varepsilon}(\partial_{x_{j}}u)\big\Vert_{L^{2}(\Omega'\times J)}\longrightarrow0}\,\,\,\,\,\,\,\,\mathit{\mathrm{as}}\,\,\,\varepsilon\to0\,,\,\,\,\,\,\,\,\,\,\,\mathrm{for\,\,each}\,\,j\in\{1,\ldots,n\}.\label{eq:Min1}
\end{equation}
$\hspace*{1em}$We now proceed to prove (\ref{eq:convergence-1}).
First observe that, for all $j\in\{1,\ldots,n\}$, we have 
\[
H_{\delta_{j}}(s)={\displaystyle \int_{0}^{s}\mathds{1}_{\{\vert w\vert\,>\,\delta_{j}\}}(\tau)\,d\tau}\,\,\,\,\,\,\,\mathrm{for\,\,every\,\,}s\in\mathbb{R}\,,
\]
\[
\sqrt{\,\widetilde{g}_{j,\varepsilon}''(\tau)}\to\mathds{1}_{\{\vert w\vert\,>\,\delta_{j}\}}(\tau)\,\,\,\,\,\,\,\mathrm{as}\,\,\varepsilon\to0\,\,\,\,\,\mathrm{for\,\,a.e.}\,\,\tau\in\mathbb{R}\,,
\]
\begin{equation}
\sqrt{\,\widetilde{g}_{j,\varepsilon}''(\tau)}\,\leq\,1\,\,\,\,\,\,\,\mathrm{for\,\,every\,\,}\tau\in\mathbb{R}\,.\label{eq:jolly2}
\end{equation}
Hence,\foreignlanguage{american}{ by the Dominated Convergence Theorem,
we conclude that, for }all $j\in\{1,\ldots,n\}$, \foreignlanguage{american}{
\[
{\displaystyle \lim_{\varepsilon\to0}\,}K_{j,\varepsilon}(s)\,=\,H_{\delta_{j}}(s)\,\,\,\,\,\,\,\,\,\,\mathrm{for\,\,every}\,\,s\in\mathbb{R}\,.
\]
}Furthermore, using (\ref{eq:H}), (\ref{eq:new_jolly}) and (\ref{eq:jolly2}),
we have,\foreignlanguage{american}{ for }all $j\in\{1,\ldots,n\}$,\begin{align*}
\left|K_{j,\varepsilon}(u_{x_{j}})-H_{\delta_{j}}(u_{x_{j}})\right|^{2}&\leq\,2\left|K_{j,\varepsilon}(u_{x_{j}})\right|^{2}+\,2\left|H_{\delta_{j}}(u_{x_{j}})\right|^{2}\\ 
&\leq\,4\,\vert u_{x_{j}}\vert^{2}\,\,\,\,\,\,\,\,\,\,\mathrm{almost\,\,everywhere\,\,in}\,\,\Omega'\times J\,.
\end{align*}Since $\vert Du\vert\in L^{2}(\Omega'\times J)$, we may apply the
Dominated Convergence Theorem again, thus obtaining 
\begin{equation}
{\displaystyle \big\Vert K_{j,\varepsilon}(u_{x_{j}})-H_{\delta_{j}}(u_{x_{j}})\big\Vert_{L^{2}(\Omega'\times J)}\longrightarrow0}\,\,\,\,\,\,\,\,\mathit{\mathrm{as}}\,\,\,\varepsilon\to0\,,\,\,\,\,\,\,\,\,\,\,\mathrm{for\,\,each}\,\,j\in\{1,\ldots,n\}.\label{eq:Min2}
\end{equation}
Finally, using Minkowski's inequality together with (\ref{eq:Min1})
and (\ref{eq:Min2}), we reach the conclusion (\ref{eq:convergence-1}).
This completes the proof.\end{proof}

\section{Proof of Theorem \ref{thm:main}\label{sec:Dimostrazione}}

\noindent $\hspace*{1em}$We are now in a position to prove Theorem
\ref{thm:main}. Indeed, using the results from Propositions \ref{prop:unif_en}
and \ref{prop:unif_en-1}, we will now show that estimate (\ref{eq:uniform})
also holds for any local weak solution $u$ of (\ref{eq:equation}),
in place of $u_{\varepsilon}$. Therefore, in the next proof we will
adopt the same assumptions and notations as in Propositions \ref{prop:unif_en}
and \ref{prop:unif_en-1}.

\noindent \begin{proof}[\bfseries{Proof of Theorem~\ref{thm:main}}]Let
$(x_{0},t_{0})\in\Omega\times I$ and $0<r<R\leq1$ be such that the
cube $Q_{R}(x_{0})$ is compactly contained in $\Omega$. In addition,
we require that $(t_{0}-R^{p},t_{0})\Subset I$. Now set 
\[
\rho\,=\,\frac{R+r}{2}
\]
and let $u_{\varepsilon}\in L^{p}(J;W^{1,p}(\Omega'))$ be the unique
weak solution to problem \textit{$(\ref{eq:regularized})$} with $\Omega'=Q_{R}(x_{0})$\foreignlanguage{american}{
and $J=(t_{0}-R^{p},t_{0})$}.\\
$\hspace*{1em}$Let us first assume that $p>2$. Then, by Proposition
\ref{prop:unif_en}, we have that 
\[
H_{\delta_{j}}(\partial_{x_{j}}u_{\varepsilon})\rightarrow H_{\delta_{j}}(\partial_{x_{j}}u)\,\,\,\,\,\,\,\,\mathit{\mathrm{strongly\,\,in}}\,\,L^{2}(Q_{R}(x_{0},t_{0}))\,\,\,\,\,\mathit{\mathrm{as}}\,\,\,\varepsilon\to0,
\]
for each $j\in\{1,\ldots,n\}$.\foreignlanguage{american}{ Thus, for
each fixed }$j\in\{1,\ldots,n\}$,\foreignlanguage{american}{ there
exists a sequence $\{\varepsilon_{k}\}_{k\in\mathbb{N}}$ such that:\vspace{0.3cm}
}

\selectlanguage{american}%
\noindent $\hspace*{1em}\bullet\,\,\,$ $0<\varepsilon_{k}<\min\,\{1,\inf\Delta^{+}\}$
for every $k\in\mathbb{N}$ and $\varepsilon_{k}\searrow0$ as $k\rightarrow+\infty$
;\vspace{0.3cm}

\noindent $\hspace*{1em}\bullet\,\,\,$ $\vert H_{\delta_{j}}(\partial_{x_{j}}u_{\varepsilon_{k}})\vert\rightarrow\vert H_{\delta_{j}}(\partial_{x_{j}}u)\vert$
almost everywhere in $Q_{R}(x_{0},t_{0})$ as $k\rightarrow+\infty$.\\

\selectlanguage{british}%
\noindent Therefore, \foreignlanguage{american}{using the definition
of $H_{\delta_{j}}$, Proposition \ref{prop:quantitative}, the fact
that }$\rho<R\leq1$\foreignlanguage{american}{, and (\ref{eq:uni})
with $\Omega'\times J=Q_{R}(x_{0},t_{0})$, we have for almost every
$z\in Q_{r}(x_{0},t_{0})$ that\begin{align*}
\vert u_{x_{j}}(z)\vert\,&\leq\,\vert H_{\delta_{j}}(u_{x_{j}}(z))\vert^{\frac{2}{p}}\,+\,\delta_{j}\,=\,\lim_{k\rightarrow\infty}\vert H_{\delta_{j}}(\partial_{x_{j}}u_{\varepsilon_{k}}(z))\vert^{\frac{2}{p}}\,+\,\delta_{j}\\
&=\,\lim_{k\rightarrow\infty}\,\max\left\{ \delta_{j},\vert\partial_{x_{j}}u_{\varepsilon_{k}}(z)\vert\right\} \\
&\leq\,\limsup_{k\rightarrow\infty}\,\,\underset{Q_{r}(x_{0},t_{0})}{\mathrm{ess\,sup}}\left(\max\left\{ \delta_{j},\vert\partial_{x_{j}}u_{\varepsilon_{k}}\vert\right\} \right)\\
&\leq\,\limsup_{k\rightarrow\infty}\,\,\max\left\{ \delta_{j},\Vert Du_{\varepsilon_{k}}\Vert_{L^{\infty}(Q_{r}(x_{0},t_{0}))}\right\} \\
&\leq\,\limsup_{k\rightarrow\infty}\,\left\{ \delta\,+\,\Vert Du_{\varepsilon_{k}}\Vert_{L^{\infty}(Q_{r}(x_{0},t_{0}))}\right\} \\
&\leq\,\limsup_{k\rightarrow\infty}\,\left\{ \delta\,+\,\frac{C}{(\rho-r)^{\vartheta p}}\left[1+\left(\iint_{Q_{\rho}(x_{0},t_{0})}\vert Du_{\varepsilon_{k}}\vert^{p}\,dx\,dt\right)^{\frac{1}{2}}\right]\right\}\\
&\leq\,\limsup_{k\rightarrow\infty}\,\,\frac{C_{1}}{(\rho-r)^{\vartheta p}}\left[1+\left(\iint_{Q_{R}(x_{0},t_{0})}\vert Du_{\varepsilon_{k}}\vert^{p}\,dx\,dt\right)^{\frac{1}{2}}\right]\\
&\leq\,\frac{C_{2}}{(R-r)^{\vartheta p}}\left[1+\left(\iint_{Q_{R}(x_{0},t_{0})}\vert Du\vert^{p}\,dx\,dt\right)^{\frac{1}{2}}\right],
\end{align*}where $\vartheta$ is defined by (\ref{eq:theta}), while the constants
$C$, $C_{1}$ and $C_{2}$ depend only on $n$, $p$ and $\delta$.
Since the above inequality holds for almost every $z\in Q_{r}(x_{0},t_{0})$,
we immediately get 
\[
\Vert u_{x_{j}}\Vert_{L^{\infty}(Q_{r}(x_{0},t_{0}))}\,\leq\,\frac{C_{2}}{(R-r)^{\vartheta p}}\left[1+\left(\iint_{Q_{R}(x_{0},t_{0})}\vert Du\vert^{p}\,dx\,dt\right)^{\frac{1}{2}}\right]<+\infty\,.
\]
This bound holds uniformly for all }$j\in\{1,\ldots,n\}$.\foreignlanguage{american}{
Thus, taking the maximum over $j$, we obtain estimate (\ref{eq:main_est}),
which proves that} $Du\in L_{loc}^{\infty}(\Omega\times I,\mathbb{R}^{n})$.\\
$\hspace*{1em}$We now consider the case $p=2$. \foreignlanguage{american}{By
(\ref{eq:convergence-1}), we conclude that, for each fixed }$j\in\{1,\ldots,n\}$,\foreignlanguage{american}{
there exists a sequence $\{\tilde{\varepsilon}_{k}\}_{k\in\mathbb{N}}$
such that:\vspace{0.3cm}
}

\selectlanguage{american}%
\noindent $\hspace*{1em}\bullet\,\,\,$ $0<\tilde{\varepsilon}_{k}<\min\,\{1,\inf\Delta^{+}\}$
for every $k\in\mathbb{N}$ and $\tilde{\varepsilon}_{k}\searrow0$
as $k\rightarrow+\infty$ ;\vspace{0.3cm}

\noindent $\hspace*{1em}\bullet\,\,\,$ $\vert K_{j,\tilde{\varepsilon}_{k}}(\partial_{x_{j}}u_{\tilde{\varepsilon}_{k}})\vert\rightarrow\vert H_{\delta_{j}}(\partial_{x_{j}}u)\vert$
almost everywhere in $Q_{R}(x_{0},t_{0})$ as $k\rightarrow+\infty$.\\

\selectlanguage{british}%
\noindent Then, \foreignlanguage{american}{arguing as above, but this
time also using the definition of $K_{j,\tilde{\varepsilon}_{k}}$
and the fact that $\sqrt{\,\widetilde{g}_{j,\tilde{\varepsilon}_{k}}''}\leq1$,
as well as (\ref{eq:uni_2}) in place of (\ref{eq:uni}), we obtain
for almost every $z\in Q_{r}(x_{0},t_{0})$ that\begin{align*}
\vert u_{x_{j}}(z)\vert\,&\leq\,\vert H_{\delta_{j}}(u_{x_{j}}(z))\vert\,+\,\delta_{j}\,=\,\lim_{k\rightarrow\infty}\vert K_{j,\tilde{\varepsilon}_{k}}(\partial_{x_{j}}u_{\tilde{\varepsilon}_{k}}(z))\vert\,+\,\delta_{j}\,\leq\,\limsup_{k\rightarrow\infty}\,\vert\partial_{x_{j}}u_{\tilde{\varepsilon}_{k}}(z)\vert\,+\,\delta\\
&\leq\,\limsup_{k\rightarrow\infty}\,\,\frac{C}{(\rho-r)^{2\vartheta}}\left[1+\left(\iint_{Q_{\rho}(x_{0},t_{0})}\vert Du_{\tilde{\varepsilon}_{k}}\vert^{2}\,dx\,dt\right)^{\frac{1}{2}}\right]+\,\delta\\
&\leq\,\limsup_{k\rightarrow\infty}\,\,\frac{C_{1}}{(R-r)^{2\vartheta}}\left[1+\left(\iint_{Q_{R}(x_{0},t_{0})}\vert Du_{\tilde{\varepsilon}_{k}}\vert^{2}\,dx\,dt\right)^{\frac{1}{2}}\right]\\
&\leq\,\frac{C_{2}}{(R-r)^{2\vartheta}}\left[1+\left(\iint_{Q_{R}(x_{0},t_{0})}\vert Du\vert^{2}\,dx\,dt\right)^{\frac{1}{2}}\right].
\end{align*}This yields the same conclusion as before and thus completes the proof.}\end{proof}

\noindent \begin{brem}A careful inspection of the proofs of Proposition
\ref{prop:quantitative} and Theorem \ref{thm:main} reveals that
the constant $C$ in (\ref{eq:main_est}) depends polynomially on
\[
\delta\,=\,1+\max\,\{\delta_{i}:i=1,\ldots,n\}\,.
\]

\noindent In particular, this constant converges to a finite quantity
depending only on $n$ and $p$ as $\max\,\{\delta_{i}\}\to0$. This
is in agreement with the fact that, in \cite{BBLV}, the constant
in the corresponding local $L^{\infty}$ estimate for $Du$ depends
only on $n$ and $p$.\end{brem}\appendix
\section{Appendix}\label{sec:appendice}$\hspace*{1em}$Let us fix $i\in\{1,\ldots,n\}$ and let $0<\varepsilon<\delta_{i}$
. We recall the $C^{2}$ function $g_{i,\varepsilon}:\mathbb{R}\to[0,\infty)$
defined in (\ref{eq:g_eps}). In this appendix, we show that $g_{i,\varepsilon}$
converges in $C^{1}(\mathbb{R})$ to 
\begin{equation}
g_{i}(s):=\,\frac{1}{2}\,(\vert s\vert-\delta_{i})_{+}^{2}\label{eq:g_i_quad}
\end{equation}
as $\varepsilon\to0$. More precisely, we establish the following
result. 
\begin{lem}
\label{lem:lemmaA1}Let $0<\varepsilon<\delta_{i}$ . Then, for every
$s\in\mathbb{R}$ we have 
\begin{equation}
\vert g_{i,\varepsilon}(s)-g_{i}(s)\vert\,\leq\,\frac{\varepsilon^{2}}{6}\,\,\,\,\,\,\,\,\,\,\,\,\,\,\mathrm{and}\,\,\,\,\,\,\,\,\,\,\,\,\,\,\vert g'_{i,\varepsilon}(s)-g'_{i}(s)\vert\,\leq\,\frac{\varepsilon}{4}\,.\label{eq:C1 convergence}
\end{equation}
\end{lem}

\noindent \begin{proof}[\bfseries{Proof}]Since $g_{i,\varepsilon}$
and $g_{i}$ are even functions, it suffices to prove the claim for
every $s\geq0$. For convenience of notation, we set 
\[
r\,=\,s-\delta_{i}\,.
\]
Recalling the definitions in (\ref{eq:g_eps}) and (\ref{eq:g_i_quad}),
we immediately have:\foreignlanguage{american}{\vspace{0.3cm}
}

\selectlanguage{american}%
\noindent $\hspace*{1em}\bullet\,\,\,$ $g_{i,\varepsilon}(r)=g_{i}(r)=0\,\,\,$
if $r\in[-\delta_{i},-\varepsilon]$;\vspace{0.3cm}

\noindent $\hspace*{1em}\bullet\,\,\,$ $\vert g_{i,\varepsilon}(r)-g_{i}(r)\vert={\displaystyle \frac{1}{12\,\varepsilon}\,(r+\varepsilon)^{3}\leq\frac{\varepsilon^{2}}{12}\,\,\,}$
if $r\in[-\varepsilon,0)$;\vspace{0.3cm}

\noindent $\hspace*{1em}\bullet\,\,\,$ $\vert g_{i,\varepsilon}(r)-g_{i}(r)\vert={\displaystyle \frac{\varepsilon^{2}}{6}\,\,\,}$
if $r\in[\varepsilon,\infty)$.\\

\selectlanguage{british}%
\noindent Furthermore, for every \foreignlanguage{american}{$r\in[0,\varepsilon]$
one has}

\selectlanguage{american}%
\noindent 
\[
\vert g_{i,\varepsilon}(r)-g_{i}(r)\vert=\left|\frac{1}{12\,\varepsilon}\,(r+\varepsilon)^{3}-\frac{1}{2}\,r^{2}\right|\leq\,\frac{\varepsilon^{2}}{6}\,,
\]
since the function 
\[
\phi(\tau):=\,\frac{1}{12\,\varepsilon}\,(\tau+\varepsilon)^{3}-\,\frac{1}{2}\,\tau^{2}
\]
is increasing and, moreover, $\phi(0)={\displaystyle \frac{\varepsilon^{2}}{12}}$
and $\phi(\varepsilon)={\displaystyle \frac{\varepsilon^{2}}{6}}$.\foreignlanguage{british}{
We have thus obtained the first inequality in (\ref{eq:C1 convergence})
for every $s\geq0$.}\\
\foreignlanguage{british}{$\hspace*{1em}$}The second inequality in\foreignlanguage{british}{
(\ref{eq:C1 convergence}) }follows by a similar argument. We leave
the details to the reader.\foreignlanguage{british}{\end{proof}}

\selectlanguage{british}%
\begin{singlespace}
\noindent \medskip{}

\noindent \textbf{Acknowledgments. }The author would like to thank
the reviewers for their valuable comments, which helped to improve
this work.\bigskip{}

\noindent \textbf{Funding.} The author is a member of the GNAMPA group
of INdAM, which partially supported his research through the INdAM--GNAMPA
2026 Project ``Esistenza e regolarità per soluzioni di equazioni
ellittiche e paraboliche anisotrope'' (CUP E53C25002010001). The
author also acknowledges financial support from the IADE\_CITTI\_2020
Project ``Intersectorial applications of differential equations''
(CUP J34I20000980006).\bigskip{}

\noindent \textbf{Data availability.} Not applicable.\addcontentsline{toc}{section}{References}
\end{singlespace}

\begin{singlespace}

\lyxaddress{\noindent \textbf{$\quad$}}
\end{singlespace}

\end{document}